\documentclass[11pt,a4paper]{article}
\usepackage[T1]{fontenc}
\usepackage{lmodern}
\usepackage[margin=1in]{geometry}
\usepackage{microtype}
\usepackage{amsmath,amssymb,amsfonts,amsthm}
\usepackage{bm}
\usepackage{multirow,booktabs,array}
\usepackage{graphicx}
\usepackage{float,placeins}
\usepackage[font=small,labelfont=bf]{caption}
\usepackage[authoryear,round]{natbib}
\usepackage{xcolor}
\usepackage{enumitem}
\usepackage{thmtools}
\definecolor{linkblue}{RGB}{20,70,130}
\usepackage[colorlinks=true,linkcolor=linkblue,citecolor=linkblue,urlcolor=linkblue]{hyperref}
\hypersetup{
  pdftitle={Synergistic Learning in Domain Adaptation with a Known Link: A Minimax Understanding},
  pdfauthor={Pengcheng Xiang, Ling Zhou, Yuhong Yang}
}
\newtheorem{theorem}{Theorem}
\newtheorem{corollary}[theorem]{Corollary}

\newtheorem{proposition}[theorem]{Proposition}

\newtheorem{definition}{Definition}
\newtheorem{remark}{Remark}

\newcommand{\mcH}{\mathcal{H}}

\newcommand{\mbE}{\mathbb{E}}

\newcommand{\Rht}{R_{h_T}}
\newcommand{\SAR}{\mathrm{SAR}(n,n_T)}

\renewcommand{\preceq}{\preccurlyeq}
\renewcommand{\succeq}{\succcurlyeq}

\title{Synergistic Learning in Domain Adaptation
with a Known Link: A Minimax Understanding}
\author{%
  Pengcheng Xiang$^{1}$ \quad Ling Zhou$^{1}$ \quad Yuhong Yang$^{2}$\\[0.6em]
  \small $^{1}$Center of Statistical Research and School of Statistics and Data Science,\\
  \small Southwest University of Finance and Economics\\[0.35em]
  \small $^{2}$Yau Math Center, Tsinghua University
}
\date{}

\begin{document}
\maketitle

\begin{abstract}
When can combining source and target data yield a much better performance
than using either one alone?
Statistically, we address this issue from a minimax viewpoint: writing $R(n,n_T)$
for the minimax risk attainable from $n$ source and $n_T$ target
observations, we ask when $R(n,n_T)$ decays \emph{strictly faster} than both
$R(n,0)$ and $R(0,n_T)$.
When this happens, we call it the \emph{Synergistic Learning Phenomenon} (SLP) and quantify
it by the \emph{Synergistic Advantage Ratio}
$\mathrm{SAR}=\min(R(n,0),R(0,n_T))/R(n,n_T)$.
In this work, we focus on domain adaptation with a \emph{known link}: the source regression function is
$f$, the target regression function is $g=\mathcal{A}f$ for a known operator
$\mathcal{A}$, and the two domains may also differ in their designs.
Three links---the identity, integration, and differentiation---together with
the behaviour of the density ratio $\rho=h_T/h_S$ determine whether
synergy is possible, at what rate, and in which sample-size window.
For the identity link (covariate shift) with a Beta$(a,1)$ source and a
uniform target, SLP is possible only if $a\geq2+1/(2\beta)$, where $\beta$
is the smoothness of $f$ (Result~I;
\citealt{zhou2022synergistic}).  An integral link $g(x)=\int_0^xf(t)\,dt$
instead can yield the joint rate $(nn_T)^{-1/2}$ even for a bounded
ratio (Result~II), and under a Beta$(a,1)$ source it has a phase diagram with
cuts at $a=2$ and $a=3$: integral SLP for $1<a<2$, an extra level-split rate
$n_T^{-1}$ for $2<a<3$, and covariate-shift-type spatial SLP for $a>3$
(Result~III).  A derivative link $g=\dot f$ synergizes only under a stronger
singularity (Result~IV), giving a sharp integration--differentiation asymmetry under a bounded ratio.
Across all three links, synergy arises from \emph{role
separation}: each dataset resolves a distinct bottleneck associated with
the link and the designs.  We characterize the sample-size windows in which
combining source and target data improves on both individual minimax rates,
and how the SAR varies with the source and target sample sizes.
Simulations support the predicted phase structure, and a mirrored
integration--differentiation analysis of real vehicle and smartphone signals
illustrates the integration--differentiation asymmetry.
\end{abstract}

\section{Introduction}
\label{sec:intro}

Can two datasets together do better than either alone? Labeled target observations are often scarce while a related source domain
supplies a far larger sample at lower cost
\citep{pan2009survey,weiss2016survey,torrey2010transfer}.  We ask whether the
two samples together can achieve a \emph{strictly better} rate of convergence
than either achieves alone---a ``$1+1>2$'' effect rather than a substitution
effect.  
In this work, we consider a domain adaptation setting where the covariate 
distribution may differ and the two regression functions are related via a known link.
We take the source regression to be $f$ and the
target to be $g=\mathcal{A}f$ for a known operator $\mathcal{A}$.
The two designs---the covariate distributions of the source and the
target---may differ with an unbounded density ratio.

Most existing rate improvements in transfer learning arise by imposing
structure on an \emph{unknown} discrepancy between the source and target
problems, for example through sparse or low-complexity parameter shifts,
function proximity, shared representations, or a smoother additive bias;
see Section~\ref{sec:related} for details.  The gain comes from a
residual that is easier to estimate than either regression function.
Our setting removes this source of gain because the known
relation $g=\mathcal{A}f$ determines $g$ from $f$, leaving no additional
regression residual to estimate.
Any ``$1+1>2$'' effect must instead come from how the known operator interacts
with the two designs: the two samples provide complementary information,
rather than an easier unknown offset.
Closest to this known-link regime are covariate shift (the identity
operator) and classical inverse problems (the same operators observed
through a single experiment); both are discussed below.

Under the identity link with a \emph{bounded} density
ratio---plain covariate shift---combining the two samples merely matches
the faster of the two individual rates \citep{pathak2022new,kpotufe2021marginal}.
An earlier work \citep{zhou2022synergistic}  showed
that  this limitation can be overcome for suitably singular designs
with an \emph{unbounded} density ratio
($\rho(x)\to\infty$ as $x\to x_0$ for some $x_0$): the combined estimator then strictly
outperforms either dataset alone in a sample-size window, the
\emph{Synergistic Learning Phenomenon} (SLP), quantified here by the
\emph{Synergistic Advantage Ratio} (SAR).

The two non-identity links have not been studied at all to the best of our knowledge.
For every $x\in[0,1]$, they take the form
\[
  g(x)=\int_0^x f(t)\,dt \quad\text{(integral link)}
  \qquad\text{or}\qquad
  g(x)=\dot f(x)=\frac{d}{dx}f(x) \quad\text{(derivative link)}.
\]
Here $f$ has H\"older smoothness $\beta$, and the goal is
to estimate $g$ under the sampling and Gaussian noise assumptions of
Section~\ref{sec:setup}.  Both experiments are informative
about~$g$, but through different channels and at different smoothness levels,
which is what makes the combined rate nonobvious.  The classical
inverse-problem literature \citep{cavalier2008nonparametric,efromovich1994adaptive}
studies the same operators with a \emph{single} dataset,
observing either the source or the target regression.

We compare the minimax risk of estimating $g$ from both
samples with the risks based on either sample alone.  
An SLP occurs when the
combined risk has a strictly faster rate than both individual risks.

As the source and target sample sizes $n$ and $n_T$ grow jointly, 
their relative scaling determines
which sample alone provides the faster minimax rate and whether
combining the samples yields a further rate improvement.
We distinguish four types of sample-size regions according to
these comparisons.
In the \emph{target dominating} (TD) region, the joint risk has the
target-only rate; in the \emph{target leading} (TL) region, the target
provides the better individual rate and gains from source data; in the
\emph{source leading} (SL) region, the source provides the better individual
rate and gains from target data; and in the \emph{source dominating} (SD)
region, the joint risk has the source-only rate.  The SLP holds exactly in
TL and SL.
Figure~\ref{fig:phase} illustrates all four regions for representative
choices of the functional link and design distributions.
In these examples, the TL and SL regions meet at the
\emph{equalizing} scale where the individual minimax rates match.
The corresponding results and Appendix~\ref{app:deferred} give the
sample-size windows and SAR formulas for the other parameter settings.
In Figure~\ref{fig:phase}(c), the joint minimax rate changes
from $n_T^{-1}$ to $(nn_T)^{-1/2}$ within the SL region, at a different
sample-size scale from that at which the individual minimax rates match.

We consider both bounded and unbounded density ratios
$\rho=h_T/h_S$ on $[0,1]$.
In the bounded-ratio setting, both design densities are bounded above
and away from zero.
In the unbounded-ratio setting, the source follows Beta$(a,1)$ with
density $h_S(x)=ax^{a-1}$, $a>1$, and the target is uniform,
so $\rho(x)=a^{-1}x^{1-a}\to\infty$ as $x\to 0$.
For these two settings, we characterize the sample-size windows,
SAR orders, and mechanisms of synergy under the identity, integral,
and derivative links.
The two main results are as follows.
\begin{itemize}[leftmargin=1.5em,topsep=2pt,itemsep=3pt]
\item \textbf{Integral link (Results~II--III).}
  Under a bounded density ratio, the joint minimax risk has rate
  $(nn_T)^{-1/2}$ in the synergistic sample-size window.
  Under a Beta$(a,1)$ source this persists for $1<a<2$; $a>2$ opens a
  second, level-split SLP at rate $n_T^{-1}$; $a>3$ replaces $(nn_T)^{-1/2}$ by
  a spatial rate; and $a=2,3$ carry logarithmic corrections
  (Theorems~\ref{thm:iii-mild}--\ref{thm:iii-spatial}).
\item \textbf{Derivative link (Result~IV).}
  Under a Beta-source design, the derivative link admits an SLP window
  with fixed polynomial sample-size exponents only when
  $a>2+\tfrac{3}{2(\beta-1)}$.
  At equality, SLP remains possible in a subpolynomial boundary layer
  (Appendix~\ref{app:deferred:iv-critical}).
\end{itemize}

These results establish an asymmetry between integration and differentiation:
in the bounded-density setting, the integral link admits synergy for
suitable relative sample sizes even when the source and target covariate
distributions coincide, whereas the derivative link does not.
Under Beta-source designs, the functional link determines the singularity
requirements, sample-size windows, and attainable SAR orders.
We characterize the joint minimax risks up to multiplicative constants
for these functional links and design distributions.
This identifies both the gains available from combining the samples and
the regimes in which no estimator can improve on the better individual
minimax rate.
We also restate the identity-link rates of \citet{zhou2022synergistic}
as Result~I, with an improved understanding at $a=2+1/(2\beta)$, so the three links can
be compared in one table.

Section~\ref{sec:related} surveys related work.
Section~\ref{sec:setup} recalls the known-link model and the SLP.
Section~\ref{sec:functional} gives Results~II--IV.
Section~\ref{sec:unified} compares the three links.
Section~\ref{sec:experiments} reports simulations and real-data
illustrations, and Section~\ref{sec:discussion} concludes.
Result~I is in Appendix~\ref{sec:covariate}; notation, additional rate and
SAR formulas, and additional experimental results are in
Appendices~\ref{app:notation}--\ref{app:semi-synthetic-unbounded}.
This preprint presents the theoretical statements and empirical results;
estimator constructions and proofs are omitted.

\section{Related Work}
\label{sec:related}

\textbf{The identity link: covariate shift.}
Under a \emph{bounded} ratio, importance-weighted ERM
\citep{cortes2010learning,mansour2009domain,sugiyama2007direct,shimodaira2000improving}
and kernel methods \citep{kpotufe2021marginal} attain the minimax rate, which
equals the faster of the two individual rates \citep{pathak2022new}: no
synergy.  Under an \emph{unbounded} ratio, \citet{zhou2022synergistic} showed
that a density singularity can produce an SLP, the first rigorous instance of
$\mathrm{SAR}\to\infty$.  Result~I restates their results for the identity
link under covariate shift.

\textbf{Non-identity links and inverse problems.}
Classical inverse-problem theory studies $g=\mathcal{A}f$ for a compact
operator \citep{cavalier2008nonparametric,engl1996regularization}, with known
single-dataset minimax rates \citep{cavalier2008nonparametric,efromovich1994adaptive}
and derivative estimation strictly harder than function estimation
\citep{rice1984bandwidth,hall1990estimation}, but observes either $f$ or $g$,
never both; whether a second dataset on the other side of a known link gives a
faster combined rate is what Results~II--IV answer.

\textbf{Unknown source--target relations.}
A broad transfer-learning literature allows the relation between source
and target problems to be unknown but restricts the discrepancy to have
additional structure.  Examples include sparse coefficient differences
\citep{li2022transfer}, aligned GLM coefficients
\citep{tian2021transfer}, shared representations
\citep{tripuraneni2021provable,du2020fewshot}, source--target function
proximity \citep{reeve2021adaptive}, close class posteriors
\citep{cai2021transfer}, small transfer bias
\citep{cai2024transfer}, and smoother additive offsets
\citep{lin2024smoothness}.  Although the structural assumptions differ,
the common mechanism is that an unknown cross-domain discrepancy is
restricted to a simpler class, making transfer statistically beneficial.
In our model, the target regression function is determined
by the source regression function through the known relation
$g=\mathcal{A}f$.  Hence any improvement from
combining the two samples must arise from how the known operator interacts
with the two designs, rather than from estimating a simpler unknown
source--target discrepancy.
Throughout we use
minimax rate analysis under H\"{o}lder smoothness
\citep{stone1982optimal,tsybakov2009introduction,donoho1990minimax,nussbaum1999equivalence},
with the SAR of Section~\ref{sec:setup} measuring when the combined rate
strictly dominates either individual rate.
\section{Problem Setup and the SLP Definition}
\label{sec:setup}
\label{sec:known-link}

We observe two independent samples, each i.i.d.\ within
its domain, with covariates in $[0,1]$ and regression models
\begin{equation}\label{eq:known-link}
  Y_i=f(X_i)+\epsilon_i,\quad X_i\sim h_S,
  \qquad
  Y_j^T=g(X_j^T)+\epsilon_j^T,\quad X_j^T\sim h_T,
  \qquad
  g=\mathcal{A}f,
\end{equation}
where $i=1,\ldots,n$ and $j=1,\ldots,n_T$.
The operator $\mathcal{A}$ is known, and the goal is to estimate $g = \mathcal{A}f$.
Throughout the minimax analysis, the noises are centered
Gaussian and independent across observations, across the two samples, and
of all covariates.  Their source and target variances are fixed, positive,
and may differ.
The design densities are fixed, with their known or
unknown status specified in each result.  For Result~II and the
bounded-ratio integration--differentiation comparison, the source density
is known; the target density need not be known, and no Lipschitz condition
on either density is imposed.  Results~III and~IV use the known
Beta-source density with known indices $(a,\beta)$ and a uniform target.
The noise assumptions above govern the two-sided minimax statements;
the weaker noise conditions for the upper bounds are stated separately
below.
We write $h$ for $h_S$ when the source design is
displayed explicitly.  The H\"older class $\mcH(\beta,\kappa)$ and the
comparison notation $\asymp,\preceq,\ll$ are collected in
Appendix~\ref{app:notation}.

For the linked class
$\mathcal F=\{(f,\mathcal Af):f\in\mcH(\beta,\kappa)\}$,
define the target $L_2(h_T)$ minimax risk
\begin{equation}\label{eq:minimax-risk}
\Rht(\mathcal F;\,n,n_T)
=\inf_{\hat g}\sup_{f\in\mcH(\beta,\kappa)}
\mbE_f\int_0^1\bigl(\hat g(x)-(\mathcal Af)(x)\bigr)^2h_T(x)\,dx .
\end{equation}
The infimum is over all estimators based on the available
samples, and $\mbE_f$ is expectation under the specified model.
We write $\Rht(g;n,n_T)$ when the linked class is clear.
Setting $n_T=0$ or $n=0$ gives the source-only risk
$\Rht(\mathcal F;n,0)$ or the target-only risk
$\Rht(\mathcal F;0,n_T)$, respectively.

\begin{definition}[Synergistic Advantage Ratio and SLP]
\label{def:slp}
The \emph{Synergistic Advantage Ratio} (SAR) at $(n,n_T)$ is
\[
  \SAR
  :=
  \frac{\min\!\bigl(\Rht(\mathcal{F};\,n,0),\;
                     \Rht(\mathcal{F};\,0,n_T)\bigr)}
       {\Rht(\mathcal{F};\,n,n_T)}.
\]
We also write $\mathrm{SAR}$ when the sample sizes are clear.
The \emph{Synergistic Learning Phenomenon} (SLP) occurs along a sample-size
sequence if $\SAR\to\infty$ as $\min(n,n_T)\to\infty$.
Since $\SAR\geq1$, $\SAR=\Theta(1)$ means no asymptotic rate improvement
over the better single-domain estimator, although constant-factor
improvements remain possible.
\end{definition}

For the identity link, Result~I records when covariate shift produces an
SLP.  The rates, including a sharper bound at $a=2+1/(2\beta)$, are in
Appendix~\ref{sec:covariate}.

\begin{figure}[t]
  \centering
  \includegraphics[width=\linewidth]{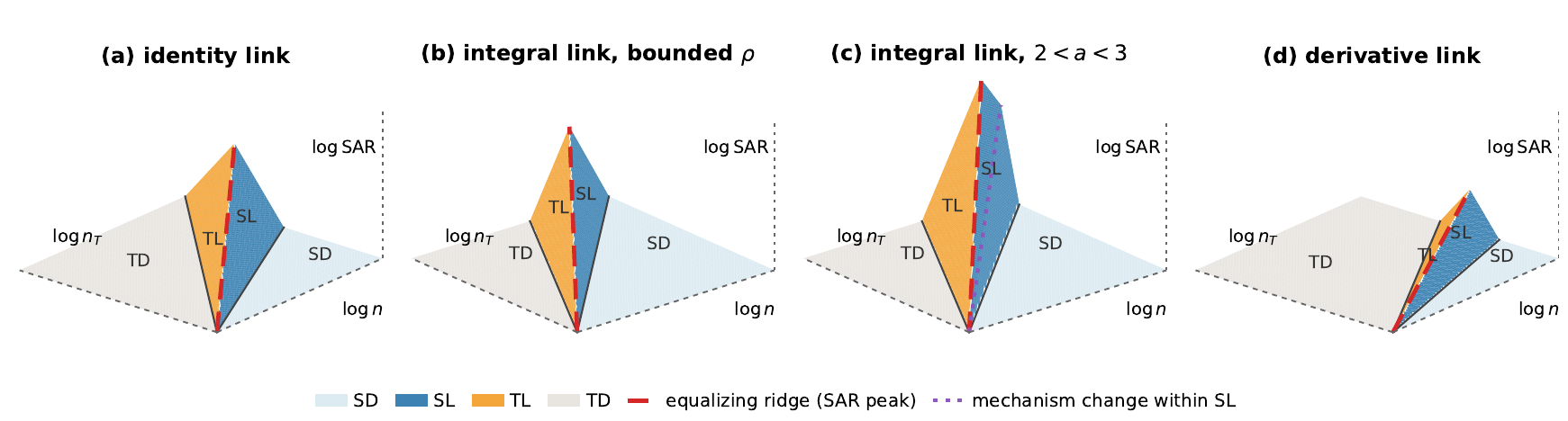}
  \caption{
Representative sample-size phase diagrams for four link--design settings:
(a) identity/Beta$(a,1)$, (b) integral/bounded density ratio,
(c) integral/Beta$(a,1)$ with $2<a<3$, and
(d) derivative/Beta$(a,1)$.
The SLP holds exactly in the SL and TL regions.
The red dashed ray marks the equalizing scale where the
individual minimax rates match.  In~(c), the purple dotted ray marks the
within-SL transition from $n_T^{-1}$ to $(nn_T)^{-1/2}$.
Parameters: $(\beta,a)=(1,4)$ in~(a), $\beta=1$ in~(b),
$(1,2.5)$ in~(c), and $(2,6)$ in~(d).
}
  \label{fig:phase}
\end{figure}

\section{Non-identity Links: Integration and Differentiation}
\label{sec:functional}

Results~II--IV below take $\mathcal{A}$ in~\eqref{eq:known-link} to be the integral
$g(x)=\int_0^xf(t)\,dt$ or the derivative $g(x)=\dot f(x)$.
The region labels in the rate displays refer to asymptotic
sample-size sequences. The SL/TL windows have strict outer endpoints to
ensure $\mathrm{SAR}\to\infty$, and share the equalizing scale where the
two individual minimax risks are of the same order.

\subsection{Result II: Integral Link with Bounded Density Ratio}
\label{sec:slp2}

\subsubsection{Setting}

The source observes $Y_i=f(X_i)+\epsilon_i$ with $X_i\sim h$ and
$f\in\mcH(\beta,\kappa)$; the target observes $Y_i^T=g(X_i^T)+\epsilon_i^T$ with
$X_i^T\sim h_T$ and $g(x)=\int_0^xf(t)\,dt\in\mcH(\beta+1,\kappa')$ for some $\kappa'>0$.

Under the conditions of Theorem~\ref{thm:ii}, the individual rates are
$\Rht(g;\,n,0)\asymp n^{-1}$ and
$\Rht(g;\,0,n_T)\asymp n_T^{-2(\beta+1)/(2\beta+3)}$.
Thus integration permits a parametric source-only rate even though the
target regression function is nonparametric.

\subsubsection{Main result}

\begin{theorem}[Minimax rate for integral link, bounded ratio]
\label{thm:ii}
Let the source and target densities satisfy
$0<c_h\leq h(x),h_T(x)\leq C_h<\infty$ for all $x\in[0,1]$,
with constants $c_h,C_h$ independent of the sample sizes.
Suppose $g(x)=\int_0^x f(t)\,dt$ and
$f\in\mcH(\beta,\kappa)$ with $\beta>0$.  Then
\begin{equation}\label{eq:ii}
  \Rht(g;\,n,n_T)
  \;\asymp\;
  \begin{cases}
    n^{-1},
      & n_T\preceq n \quad\text{(SD)},\\[4pt]
    (n\,n_T)^{-1/2},
      & n\ll n_T\preceq n^{\frac{2\beta+3}{2(\beta+1)}}
        \quad\text{(SL)},\\[4pt]
    (n\,n_T)^{-1/2},
      & n^{\frac{2\beta+3}{2(\beta+1)}}\preceq n_T
        \ll n^{\frac{2\beta+3}{2\beta+1}} \quad\text{(TL)},\\[4pt]
    n_T^{-\frac{2(\beta+1)}{2\beta+3}},
      & n_T\succeq n^{\frac{2\beta+3}{2\beta+1}}
        \quad\text{(TD)}.
  \end{cases}
\end{equation}
\end{theorem}

Here the source density $h$ is assumed known; the target density $h_T$
need not be.
The ratio being upper bounded, there is no covariate-shift singularity and
the SLP has a \emph{different} mechanism.

\begin{corollary}[SAR for integral link, bounded ratio]
\label{cor:sar-ii}
Under the conditions of Theorem~\ref{thm:ii}, the SLP occurs
if and only if $n\ll n_T\ll n^{(2\beta+3)/(2\beta+1)}$.
Within this window,
\begin{equation}\label{eq:sar-ii}
  \SAR\asymp
  \begin{cases}
    (n_T n^{-1})^{1/2},
    & n\ll n_T\preceq n^{\frac{2\beta+3}{2(\beta+1)}}
      \quad\text{(SL)},\\[4pt]
    \bigl(n\,n_T^{-\frac{2\beta+1}{2\beta+3}}\bigr)^{1/2},
    & n^{\frac{2\beta+3}{2(\beta+1)}}\preceq n_T
      \ll n^{\frac{2\beta+3}{2\beta+1}} \quad\text{(TL)}.
  \end{cases}
\end{equation}
In the SD and TD regions, $\mathrm{SAR}=\Theta(1)$.
\end{corollary}

The source-only risk already has the parametric
rate $n^{-1}$, so synergy requires $n\ll n_T$.
Within the SLP window, the SL and TL regions meet at
$n_T\asymp n^{(2\beta+3)/(2(\beta+1))}$, where the SAR attains its
largest order.

\subsection{Result III: Integral Link with Unbounded Density Ratio}
\label{sec:slp3}

\begingroup
\subsubsection{Setting}
\label{sec:iii-setting}
\label{app:deferred:iii-basis}

For $a>1$ the source observes $Y_i=f(X_i)+\epsilon_i$ with $X_i\sim h(x)=ax^{a-1}$
(Beta$(a,1)$), $f\in\mcH(\beta,\kappa)$ and anchor $g(0)=0$, while the target
observes $Y_i^T=g(X_i^T)+\epsilon_i^T$ with $X_i^T\sim h_T$, $h_T\equiv1$, and
$g(x)=\int_0^xf(t)\,dt$: the integral link of Result~II combined with the
singular design of Result~I.
The source density $h$ and the indices $(a,\beta)$ are assumed known,
with $\beta>0$.

The Result~III rate statements are unconditional minimax equivalences
for every fixed $a>1$ under the model and regularity conditions of
this setting.

\begin{remark}[Relaxed density assumption for Result~III]
\label{rem:iii-density}
Result~III uses $h(x)=ax^{a-1}$ with $a>1$, so $h\to0$ and
$h_T/h\to\infty$ at the origin.  The positive lower bound on the source
density in Theorem~\ref{thm:ii} therefore fails.
The range $a\leq1$, where this density-ratio singularity is absent, is not
part of Result~III.
\end{remark}

\subsubsection{Mild singularity: integral SLP only}
\label{sec:iii-mild}

\begin{theorem}[Integral SLP for $1<a<2$]
\label{thm:iii-mild}
\label{thm:iii}
Let $h$ and $h_T$ be as defined above, and let $1<a<2$.
Suppose $g(x)=\int_0^x f(t)\,dt$, $f\in\mcH(\beta,\kappa)$, $\beta>0$, and
$g(0)=0$.
Then
\begin{equation}\label{eq:iii-mild}
  \Rht(g;\,n,n_T)\;\asymp\;
  \begin{cases}
    n^{-1}, & n_T\preceq n \quad\text{(SD)},\\[4pt]
    (nn_T)^{-1/2},
      & n\ll n_T\preceq n^{\frac{2\beta+3}{2(\beta+1)}}
        \quad\text{(SL)},\\[4pt]
    (nn_T)^{-1/2},
      & n^{\frac{2\beta+3}{2(\beta+1)}}\preceq n_T
        \ll n^{\frac{2\beta+3}{2\beta+1}} \quad\text{(TL)},\\[4pt]
    n_T^{-\frac{2(\beta+1)}{2\beta+3}},
      & n_T\succeq n^{\frac{2\beta+3}{2\beta+1}}
        \quad\text{(TD)}.
  \end{cases}
\end{equation}
\end{theorem}

Here the source density $h$ and the indices $(a,\beta)$ are assumed known.
Although $\rho\to\infty$, \eqref{eq:iii-mild} coincides
with the rate in Result~II. The only synergistic rate is $(nn_T)^{-1/2}$,
on $n\ll n_T\ll n^{(2\beta+3)/(2\beta+1)}$.

\subsubsection{Intermediate singularity: level-split + integral SLP}
\label{sec:iii-level}
\label{app:deferred:iii-mid}

\begin{theorem}[Extra level-split SLP for $2<a<3$]
\label{thm:iii-level}
Let $h$ and $h_T$ be as defined above, and let $2<a<3$.
Suppose $g(x)=\int_0^x f(t)\,dt$, $f\in\mcH(\beta,\kappa)$, $\beta>0$, and
$g(0)=0$.
Then
\begin{equation}\label{eq:iii-level}
  \Rht(g;\,n,n_T)\;\asymp\;
  \begin{cases}
    n^{-\frac{2(\beta+1)}{2\beta+a}},
      & n_T\preceq n^{\frac{2(\beta+1)}{2\beta+a}}
        \quad\text{(SD)},\\[4pt]
    n_T^{-1},
      & n^{\frac{2(\beta+1)}{2\beta+a}}\ll n_T\preceq n
        \quad\text{(SL)},\\[4pt]
    (nn_T)^{-1/2},
      & n\preceq n_T\preceq n^{\frac{2\beta+3}{2\beta+a}}
        \quad\text{(SL)},\\[4pt]
    (nn_T)^{-1/2},
      & n^{\frac{2\beta+3}{2\beta+a}}\preceq n_T
        \ll n^{\frac{2\beta+3}{2\beta+1}} \quad\text{(TL)},\\[4pt]
    n_T^{-\frac{2(\beta+1)}{2\beta+3}},
      & n_T\succeq n^{\frac{2\beta+3}{2\beta+1}}
        \quad\text{(TD)}.
  \end{cases}
\end{equation}
\end{theorem}

\begin{remark}[Why the source-only rate is slower than $n^{-1}$ for $a>2$]
\label{rem:iii-no-ninv}
Source data observe $f$, while the target parameter is the running integral
$g(x)=\int_0^x f(t)\,dt$.
For $a>2$, the sparse source design near the origin prevents parametric
source-only estimation of this integral functional.  In this range,
$
  \Rht(g;\,n,0)
  \;\succeq\;
  n^{-2(\beta+1)/(2\beta+a)}.
$
For $a>2$ the right-hand side is strictly slower than~$n^{-1}$.
For $1<a<2$, the source-only rate remains $n^{-1}$, as recorded in
Theorem~\ref{thm:iii-mild}.
\end{remark}

\begin{remark}[Role split and the two synergistic rates]
\label{rem:iii-nt-inv}
Write the risk as the sum of two pieces: the running-integral functionals
of~$g$ that are irregular for the source when $a>2$, and the remaining
structure of~$g$.

\begin{itemize}[leftmargin=*,topsep=2pt,itemsep=2pt]
\item \textbf{Irregular running-integral functionals.}
  By Remark~\ref{rem:iii-no-ninv}, source-only estimation of these
  functionals has rate $n^{-2(\beta+1)/(2\beta+a)}$.
  The target observes~$g$ directly, so the same functionals are ordinary
  parametric problems with rate~$n_T^{-1}$.
  Target alone estimating the \emph{whole} of
  $g\in\mcH(\beta+1)$ would still cost the nonparametric rate
  $n_T^{-2(\beta+1)/(2\beta+3)}$; the rate~$n_T^{-1}$ applies only to
  this finite-dimensional piece.

\item \textbf{Remaining structure.}
  For $a<3$ this part remains parametrically estimable from the source
  (observations of~$f$), at rate~$n^{-1}$.
  The Result~II-type combination with the target improves it further to
  $(nn_T)^{-1/2}$, but only when $n_T\gtrsim n$; for $n_T\preceq n$ one has
  $(nn_T)^{-1/2}\succeq n^{-1}$, so the source-only rate~$n^{-1}$ is
  already optimal for this piece.
\end{itemize}

\end{remark}

\subsubsection{Strong singularity: level-split + covariate-shift SLP}
\label{sec:iii-spatial}
\label{app:deferred:iii-strong}

\begin{theorem}[Spatial singularity SLP for $a>3$]
\label{thm:iii-spatial}
Under the setting of Section~\ref{sec:iii-setting} with $a>3$,
\begin{equation}\label{eq:iii-spatial}
  \Rht(g;\,n,n_T)\;\asymp\;
  \begin{cases}
    n^{-\frac{2(\beta+1)}{2\beta+a}},
      & n_T\preceq n^{\frac{2(\beta+1)}{2\beta+a}}
        \quad\text{(SD)},\\[4pt]
    n_T^{-1},
      & n^{\frac{2(\beta+1)}{2\beta+a}}\ll n_T
        \preceq n^{\frac{2\beta+3}{2\beta+a}} \quad\text{(SL)},\\[4pt]
    \Bigl(n_T^{\frac{2\beta+1}{2\beta+3}}n^{-1}\Bigr)^{\frac{1}{a-1}}
    n_T^{-\frac{2(\beta+1)}{2\beta+3}},
      & n^{\frac{2\beta+3}{2\beta+a}}\preceq n_T
        \ll n^{\frac{2\beta+3}{2\beta+1}} \quad\text{(TL)},\\[4pt]
    n_T^{-\frac{2(\beta+1)}{2\beta+3}},
      & n_T\succeq n^{\frac{2\beta+3}{2\beta+1}}
        \quad\text{(TD)}.
  \end{cases}
\end{equation}
\end{theorem}

\begin{remark}[Why the mid-rate is covariate-shift type, and why
$(nn_T)^{-1/2}$ disappears]
\label{rem:iii-spatial}
For $a>3$, even after the irregular running-integral functionals are
handled as in Remark~\ref{rem:iii-nt-inv}, the source still cannot
estimate arbitrarily fine local features of~$g$ near~$x=0$ from
observations of~$f$.
Consequently there is no Result~II-type regime in which the remaining
structure of~$g$ can be combined at the integral rate~$(nn_T)^{-1/2}$.
The two samples instead contribute in different spatial regions, as in
Result~I, with source
information for levels of~$g$ entering through the derivative link
$f=\dot g$.
The resulting mid-rate in~\eqref{eq:iii-spatial} is therefore a
covariate-shift-type synergistic rate, not an integral-type rate.

By contrast, in Theorem~\ref{thm:iii-level} ($2<a<3$), those local features
remain parametrically estimable from the source.
The integral combination $(nn_T)^{-1/2}$ is then available, and no
covariate-shift-type mid-rate appears.
The cut $a=3$ is exactly the transition between these two behaviours.
\end{remark}

\subsubsection{SAR for \texorpdfstring{$1<a<2$, $2<a<3$, and $a>3$}{1<a<2, 2<a<3, and a>3}}
\label{sec:iii-sar}

\begin{corollary}[SAR for Result~III when $a\notin\{2,3\}$]
\label{cor:sar-iii}
Under the setting of Section~\ref{sec:iii-setting}, suppose
$a\in(1,2)\cup(2,3)\cup(3,\infty)$.
Table~\ref{tab:iii-phase} summarizes the SLP windows and synergistic rates.
On each such window $\mathrm{SAR}\to\infty$ as below;
in the SD and TD regions, $\mathrm{SAR}=\Theta(1)$.
\begin{subequations}\label{eq:sar-iii}
For $1<a<2$, the SAR coincides with~\eqref{eq:sar-ii}:
\begin{equation}\label{eq:sar-iii-mild}
  \SAR\asymp
  \begin{cases}
    (n_T n^{-1})^{1/2},
      & n\ll n_T\preceq n^{\frac{2\beta+3}{2(\beta+1)}}
        \quad\text{(SL)},\\[4pt]
    \bigl(n\,n_T^{-\frac{2\beta+1}{2\beta+3}}\bigr)^{1/2},
      & n^{\frac{2\beta+3}{2(\beta+1)}}\preceq n_T
        \ll n^{\frac{2\beta+3}{2\beta+1}} \quad\text{(TL)}.
  \end{cases}
\end{equation}
For $2<a<3$,
\begin{equation}\label{eq:sar-iii-level}
  \SAR\asymp
  \begin{cases}
    n_T n^{-\frac{2(\beta+1)}{2\beta+a}},
      & n^{\frac{2(\beta+1)}{2\beta+a}}\ll n_T\preceq n
        \quad\text{(SL)},\\[4pt]
    n^{-\frac{2(\beta+1)}{2\beta+a}+\frac12}\,n_T^{1/2},
      & n\preceq n_T\preceq n^{\frac{2\beta+3}{2\beta+a}}
        \quad\text{(SL)},\\[4pt]
    \bigl(n\,n_T^{-\frac{2\beta+1}{2\beta+3}}\bigr)^{1/2},
      & n^{\frac{2\beta+3}{2\beta+a}}\preceq n_T
        \ll n^{\frac{2\beta+3}{2\beta+1}} \quad\text{(TL)}.
  \end{cases}
\end{equation}
For $a>3$,
\begin{equation}\label{eq:sar-iii-spatial}
  \SAR\asymp
  \begin{cases}
    n_T n^{-\frac{2(\beta+1)}{2\beta+a}},
      & n^{\frac{2(\beta+1)}{2\beta+a}}\ll n_T
        \preceq n^{\frac{2\beta+3}{2\beta+a}} \quad\text{(SL)},\\[4pt]
    \bigl(n\,n_T^{-\frac{2\beta+1}{2\beta+3}}\bigr)^{\frac{1}{a-1}},
      & n^{\frac{2\beta+3}{2\beta+a}}\preceq n_T
        \ll n^{\frac{2\beta+3}{2\beta+1}} \quad\text{(TL)}.
  \end{cases}
\end{equation}
\end{subequations}
In~\eqref{eq:sar-iii-level}, the change of rate at $n_T\asymp n$
lies within SL.  In~\eqref{eq:sar-iii-spatial}, the spatial rate begins
at the equalizing scale $n_T\asymp n^{(2\beta+3)/(2\beta+a)}$;
the target-only risk is therefore the smaller individual benchmark
throughout TL, up to constants at that scale.
\end{corollary}

\subsubsection{The cases \texorpdfstring{$a=2$ and $a=3$}{a=2 and a=3}, phase diagram, and SAR}
\label{sec:iii-critical}
\label{app:deferred:iii-critical}
\label{sec:iii-phase}

\begin{proposition}[Minimax rate at $a=2$]
\label{prop:iii-a2}
At $a=2$,
\begin{equation}\label{eq:iii-a2}
  \Rht(g;\,n,n_T)
  \;\asymp\;
  \begin{cases}
    (\log n)\,n^{-1},
      & n_T\preceq n(\log n)^{-1} \quad\text{(SD)},\\[4pt]
    n_T^{-1},
      & n(\log n)^{-1}\ll n_T\preceq n \quad\text{(SL)},\\[4pt]
    (nn_T)^{-1/2},
      & n\preceq n_T\preceq
        \bigl(n(\log n)^{-1}\bigr)^{\frac{2\beta+3}{2(\beta+1)}}
        \quad\text{(SL)},\\[4pt]
    (nn_T)^{-1/2},
      & \bigl(n(\log n)^{-1}\bigr)^{\frac{2\beta+3}{2(\beta+1)}}
        \preceq n_T\ll n^{\frac{2\beta+3}{2\beta+1}}
        \quad\text{(TL)},\\[4pt]
    n_T^{-\frac{2(\beta+1)}{2\beta+3}},
      & n_T\succeq n^{\frac{2\beta+3}{2\beta+1}}
        \quad\text{(TD)}.
  \end{cases}
\end{equation}
\end{proposition}

\noindent
Relative to the limits $a\to 2^\pm$, the leading power remains $n^{-1}$,
an extra $\log n$ appears in the source-dominated risk, and the
$n_T^{-1}$ window shrinks to $n/\log n\ll n_T\preceq n$.
SLP holds on $n/\log n\ll n_T\ll n^{(2\beta+3)/(2\beta+1)}$;
SL and TL meet at $n_T\asymp(n/\log n)^{(2\beta+3)/(2(\beta+1))}$.

\begin{proposition}[Minimax rate at $a=3$]
\label{prop:iii-a3}
At $a=3$,
\begin{equation}\label{eq:iii-a3}
  \Rht(g;\,n,n_T)
  \;\asymp\;
  \begin{cases}
    n^{-\frac{2(\beta+1)}{2\beta+3}},
      & n_T\preceq n^{\frac{2(\beta+1)}{2\beta+3}}
        \quad\text{(SD)},\\[4pt]
    n_T^{-1},
      & n^{\frac{2(\beta+1)}{2\beta+3}}\ll n_T
        \preceq n(\log n)^{-1} \quad\text{(SL)},\\[4pt]
    (\log n)\,n^{-1},
      & n(\log n)^{-1}\preceq n_T\preceq n \quad\text{(SL)},\\[4pt]
    (nn_T)^{-1/2}\Bigl(1+\log_+\bigl(en\,
      n_T^{-\frac{2\beta+1}{2\beta+3}}\bigr)\Bigr),
      & n\preceq n_T\ll n^{\frac{2\beta+3}{2\beta+1}}
        \quad\text{(TL)},\\[4pt]
    n_T^{-\frac{2(\beta+1)}{2\beta+3}},
      & n_T\succeq n^{\frac{2\beta+3}{2\beta+1}}
        \quad\text{(TD)}.
  \end{cases}
\end{equation}
\end{proposition}

\noindent
The positive-part logarithm keeps the factor bounded away
from zero when $n_T\asymp n^{(2\beta+3)/(2\beta+1)}$, where the last two rate
expressions agree up to constants.
The rate $(nn_T)^{-1/2}$ survives only with a logarithmic factor;
the polynomial spatial mid-rate of Theorem~\ref{thm:iii-spatial}
requires $a>3$.
SLP holds on
$n^{2(\beta+1)/(2\beta+3)}\ll n_T\ll n^{(2\beta+3)/(2\beta+1)}$.
The transition at $n_T\asymp n/\log n$ lies within SL, while the
equalizing scale $n_T\asymp n$ separates SL and TL.

\begin{remark}[Left and right of the cuts $a=2$ and $a=3$]
\label{rem:iii-cuts}
Theorems~\ref{thm:iii-mild}--\ref{thm:iii-spatial} are organised by how
well the \emph{source} sample alone can estimate~$g$.
\begin{enumerate}[leftmargin=*,topsep=2pt,itemsep=2pt]
\item \textbf{Cut $a=2$.}
  As explained in Remark~\ref{rem:iii-no-ninv}, when $a>2$ the source
  design puts too little mass near the origin for the integral map
  $f\mapsto g$ to remain parametrically estimable, so
  $\Rht(g;n,0)$ is of order $n^{-2(\beta+1)/(2\beta+a)}$ rather
  than~$n^{-1}$.
  The target sample observes~$g$ directly and attains risk~$n_T^{-1}$ for
  the corresponding finite-dimensional functionals of~$g$.
  This produces the extra synergistic window of rate $n_T^{-1}$ in
  Theorem~\ref{thm:iii-level}.

\item \textbf{Cut $a=3$.}
  Independently of the running-integral functionals above, the source must
  still estimate local features of~$g$ near~$x=0$ from observations of
  $f=\dot g$.
  If $a<3$, arbitrarily fine features remain estimable at parametric
  source rates and the
  Result~II-type combination $(nn_T)^{-1/2}$ continues to apply.
  If $a>3$, the source lacks this local information near the
  singularity; the samples contribute in different spatial regions,
  with the mid-rate of
  Theorem~\ref{thm:iii-spatial} (the same geometry as Result~I, but for
  the integral link).
\end{enumerate}
Thus unboundedness of $\rho=h_T/h$ (which already holds for all $a>1$)
is not by itself enough to change the SLP structure: the first change
occurs when source-only estimation of~$g$ is forced below the parametric
rate~$n^{-1}$ ($a=2$), and the second when local features of~$g$ near~$0$
cease to be parametrically estimable from the source ($a=3$).
\end{remark}

\begin{remark}[Why the rates change exactly at $a=2$ and $a=3$]
\label{rem:iii-critical-rates}
Propositions~\ref{prop:iii-a2}--\ref{prop:iii-a3} record the
rates at $a=2$ and $a=3$.  They differ from both neighbouring regimes
because each cut is a borderline second-moment condition.
At $a=2$, the source-dominated rate is $(\log n)\,n^{-1}$ and the
level-split window starts at $n_T\asymp n/\log n$.
At $a=3$, the joint risk has a band of rate $(\log n)\,n^{-1}$,
followed by a log-corrected integral rate $(nn_T)^{-1/2}$.
The polynomial spatial mid-rate appears only for $a>3$.
These logarithms distinguish the cases $a=2$ and $a=3$ from the
neighbouring parameter ranges without changing their leading powers.
\end{remark}

The SAR formulas for $1<a<2$, $2<a<3$, and $a>3$ are given in
Corollary~\ref{cor:sar-iii}.
The SAR formulas at $a=2,3$ remain in
Corollary~\ref{cor:sar-iii-critical} (Appendix~\ref{app:deferred:iii-sar}),
and Table~\ref{tab:iii-phase} (Appendix~\ref{app:deferred:iii-table})
collects the SLP windows and synergistic rates for every range of $a$.
In the SD and TD regions, $\mathrm{SAR}=\Theta(1)$.

\endgroup

\subsection{Result IV: Derivative Link with Unbounded Density Ratio}
\label{sec:slp4}

\subsubsection{Setting}

The source observes $Y_i=f(X_i)+\epsilon_i$ with $X_i\sim h(x)=ax^{a-1}$, $a>1$, and
$f\in\mcH(\beta,\kappa)$, $\beta>1$; the target observes
$Y_i^T=\dot f(X_i^T)+\epsilon_i^T$ with $X_i^T\sim h_T$, $h_T\equiv1$, and
$\dot f\in\mcH(\beta-1,\kappa'')$, so differentiation lowers smoothness by~1
while the design carries the Result~I singularity.

The target-only rate is
$\Rht(\dot f;0,n_T)\asymp n_T^{-2(\beta-1)/(2\beta-1)}$, while the source-only rate is
$\Rht(\dot f;n,0)\asymp n^{-2(\beta-1)/(2\beta+1)}$ for
$a<2+\tfrac{3}{2(\beta-1)}$, carries a $\log n$ factor at equality, and slows to
$n^{-(2\beta-1)/(2\beta+a)}$ above it.

\subsubsection{Main result}

\begin{theorem}[Minimax rate for derivative link under Beta$(a,1)$]
\label{thm:iv}
Let $a>1$, $h(x)=ax^{a-1}$, $h_T\equiv1$, and
$f\in\mcH(\beta,\kappa)$, $\beta>1$.  The following rates hold.
\begin{itemize}[leftmargin=*,topsep=2pt,itemsep=2pt,parsep=0pt]
\item If $a<2+\tfrac{3}{2(\beta-1)}$:
  $\Rht(\dot{f};\,n,n_T)\asymp
  \min\bigl(n^{-\frac{2(\beta-1)}{2\beta+1}},\,
            n_T^{-\frac{2(\beta-1)}{2\beta-1}}\bigr)$.
\item If $a>2+\tfrac{3}{2(\beta-1)}$:
      \begin{equation}\label{eq:iv}
        \Rht(\dot{f};\,n,n_T)\;\asymp\;
        \begin{cases}
          n^{-\frac{2\beta-1}{2\beta+a}},
            & n_T\preceq n^{\frac{2\beta-1}{2\beta+a}}
              \quad\text{(SD)},\\[4pt]
          \Bigl(n_T^{\frac{2\beta+1}{2\beta-1}}n^{-1}\Bigr)^{\frac{1}{a-1}}
          n_T^{-\frac{2(\beta-1)}{2\beta-1}},
            & n^{\frac{2\beta-1}{2\beta+a}}\ll n_T
              \preceq n^{\frac{(2\beta-1)^2}{2(\beta-1)(2\beta+a)}}
              \quad\text{(SL)},\\[4pt]
          \Bigl(n_T^{\frac{2\beta+1}{2\beta-1}}n^{-1}\Bigr)^{\frac{1}{a-1}}
          n_T^{-\frac{2(\beta-1)}{2\beta-1}},
            & n^{\frac{(2\beta-1)^2}{2(\beta-1)(2\beta+a)}}\preceq n_T
              \ll n^{\frac{2\beta-1}{2\beta+1}} \quad\text{(TL)},\\[4pt]
          n_T^{-\frac{2(\beta-1)}{2\beta-1}},
            & n_T\succeq n^{\frac{2\beta-1}{2\beta+1}}
              \quad\text{(TD)}.
        \end{cases}
      \end{equation}
\end{itemize}
\end{theorem}

The indices $(a,\beta)$ are assumed known.
The source observes function levels while the target observes derivatives,
so the two samples provide information about the target at different
smoothness levels.
If $a<2+\tfrac{3}{2(\beta-1)}$, there is no singularity-driven SLP
(same obstruction as Theorem~\ref{thm:deriv}), and no parametric IPW
fallback: differentiation admits no exact-unbiasedness channel.
If $a>2+\tfrac{3}{2(\beta-1)}$, the rates are \emph{not} Result~I with
$\beta\mapsto\beta-1$. On a window of scale $t$ and source probability
mass $H$, source levels contribute Fisher information of order $nHt^{2}$
rather than $nH$.
The case $a=2+\tfrac{3}{2(\beta-1)}$ carries a logarithmic
correction and is stated in Appendix~\ref{app:deferred:iv-critical2}.
The source-only rate under the Beta design, and why the SLP and SAR conclusions
are unaffected by it, are recorded in Appendix~\ref{app:deferred:iv-src}.

\begin{corollary}[SAR for derivative link under Beta$(a,1)$]
\label{cor:sar-iv}
Under the conditions of Theorem~\ref{thm:iv},
for $a>2+\tfrac{3}{2(\beta-1)}$ the SLP occurs iff
$n^{(2\beta-1)/(2\beta+a)}\ll n_T\ll n^{(2\beta-1)/(2\beta+1)}$.
Within this window,
\begin{equation}\label{eq:sar-iv}
  \SAR\asymp
  \begin{cases}
    \bigl(n_T\,n^{-\frac{2\beta-1}{2\beta+a}}\bigr)^{
      \frac{2(\beta-1)}{2\beta-1}
      -\frac{2\beta+1}{(2\beta-1)(a-1)}},
    & n^{\frac{2\beta-1}{2\beta+a}}\ll n_T
      \preceq n^{\frac{(2\beta-1)^2}{2(\beta-1)(2\beta+a)}}
      \quad\text{(SL)},\\[4pt]
    \bigl(n\,n_T^{-\frac{2\beta+1}{2\beta-1}}\bigr)^{\frac{1}{a-1}},
    & n^{\frac{(2\beta-1)^2}{2(\beta-1)(2\beta+a)}}\preceq n_T
      \ll n^{\frac{2\beta-1}{2\beta+1}} \quad\text{(TL)}.
  \end{cases}
\end{equation}
The SAR attains its largest order at the SL--TL equalizing scale
$n_T\asymp n^{(2\beta-1)^2/(2(\beta-1)(2\beta+a))}$.
In the SD and TD regions, $\mathrm{SAR}=\Theta(1)$.
For $a<2+\tfrac{3}{2(\beta-1)}$, $\SAR=\Theta(1)$.
At $a=2+\tfrac{3}{2(\beta-1)}$, the SLP occurs precisely
in the subpolynomial boundary layer specified in
Appendix~\ref{app:deferred:iv-critical}.
\end{corollary}

Compared with Results~I and~III, the threshold
$a>2+\tfrac{3}{2(\beta-1)}$ is the \emph{strictest}: estimating $\dot{f}$
from source levels incurs the extra Fisher factor $t^{2}$, so a stronger
singularity is needed than for a common $\mcH(\beta)$ regression.
At equality the SAR is only subpolynomially large, so there is no fixed
polynomial-exponent SLP window; see Appendix~\ref{app:deferred:iv-critical}.

\subsection{Asymmetry Under Bounded Ratio}
\label{sec:asymmetry}

Results~II--IV show that the density ratio and the functional link between
the source and target regression functions jointly determine whether SLP is possible.
To isolate the effect of the functional link, we now restrict to the bounded-density
setting of Result~II and compare the integral and derivative links.
Under these conditions, the integral link admits SLP for suitable relative
sample sizes (Result~II) but the derivative link does not---a sharp asymmetry
between the two links.

For the derivative link, the source observes
$Y_i=f(X_i)+\epsilon_i$ with $X_i\sim h$ and $f\in\mcH(\beta,\kappa)$, $\beta>1$,
and the target observes $Y_i^T=\dot f(X_i^T)+\epsilon_i^T$ with $X_i^T\sim h_T$ and
$\dot f\in\mcH(\beta-1,\kappa'')$.

\begin{theorem}[No SLP for derivative link under bounded ratio]
\label{thm:deriv}
Let the source and target densities satisfy
$0<c_h\leq h(x),h_T(x)\leq C_h<\infty$ for all $x\in[0,1]$,
with constants $c_h,C_h$ independent of the sample sizes.
Let $\beta>1$ and use the linked class
$\mathcal F_{\beta,\kappa}^{\mathrm{BD}}
=\{(f,g):f\in\mcH(\beta,\kappa),\ g=f'\}$.
Write $\mathcal R_{\mathrm{BD}}^{\mathrm G}(n,n_T)
:=\Rht(\mathcal F_{\beta,\kappa}^{\mathrm{BD}};n,n_T)$ for its minimax risk
under the Gaussian model in Section~\ref{sec:setup}.  Then
\[
  \mathcal R_{\mathrm{BD}}^{\mathrm G}(n,n_T)
  \;\asymp\;
  \min\!\bigl(n^{-\frac{2(\beta-1)}{2\beta+1}},\;
             n_T^{-\frac{2(\beta-1)}{2\beta-1}}\bigr).
\]
\end{theorem}

$\mathcal R_{\mathrm{BD}}^{\mathrm G}(n,n_T)$ is the minimum of the two
individual risks up to constants, so $\SAR=\Theta(1)$ for every sample-size
sequence, and the SLP does not occur.

The upper bounds for Results~II--IV and
Theorem~\ref{thm:deriv}, in all stated parameter cases, remain valid for
independent observations with conditionally mean-zero noises and uniformly
bounded conditional second moments.

\begin{remark}[Integration--differentiation asymmetry under bounded ratio]
\label{rem:asymmetry}
\label{thm:asymmetry}
Comparing Theorem~\ref{thm:ii} with Theorem~\ref{thm:deriv} gives a sharp
integration--differentiation asymmetry under their common density conditions, with the source density known
in both experiments: for $g=\int_0^\cdot f$ the SLP occurs on the TL/SL window
of Corollary~\ref{cor:sar-ii} and $\SAR\to\infty$, whereas for $g=\dot f$ it
occurs for no $(n,n_T)$ and $\SAR=\Theta(1)$.
The barrier is specific to the bounded-ratio regime.
Under the Beta model, the integral link admits SLP for every $a>1$
(Result~III).  The derivative link requires $a>2+\tfrac{3}{2(\beta-1)}$
for a polynomial SLP window, with a subpolynomial SLP window at equality
(Result~IV).
\end{remark}

Integration raises smoothness and admits an exactly unbiased estimator, whereas
differentiation lowers smoothness and amplifies source noise
(Appendix~\ref{app:deferred:asym}).

\section{Unified Characterization of the SLP}
\label{sec:unified}

As summarized in Table~\ref{tab:summary}, the identity link needs a
strong singularity, the integral link none, and the derivative link a strictly
stronger one. The three functional links therefore have different
singularity requirements for SLP.

\begin{table}[H]
\centering
\caption{Summary of Results~I--IV and the bounded-ratio negative result;
$\mathfrak C_\beta(n,n_T)$ is the identity-link rate at $a=2+1/(2\beta)$
in~\eqref{eq:cs-critical}.}
\label{tab:summary}
\scriptsize
\setlength{\tabcolsep}{3pt}
\renewcommand{\arraystretch}{1.05}
\begin{tabular}{@{}>{\raggedright\arraybackslash}p{0.055\textwidth}
                    >{\raggedright\arraybackslash}p{0.17\textwidth}
                    >{\raggedright\arraybackslash}p{0.27\textwidth}
                    >{\raggedright\arraybackslash}p{0.27\textwidth}
                    >{\raggedright\arraybackslash}p{0.13\textwidth}@{}}
\toprule
Result & Setting & Condition for SLP & Synergistic rate & SAR \\
\midrule
I & Identity link (covariate shift), $\rho\to\infty$ & $a>2+\frac{1}{2\beta}$:
  polynomial SL/TL window;
  equality: $\log(n/n_T)=o(\log n)$ and $n_T/n\to0$
  & above threshold:
  $\bigl(n_T n^{-1}\bigr)^{\!\frac{1}{a-1}}n_T^{-\frac{2\beta}{2\beta+1}}$;
  equality: $\mathfrak C_\beta(n,n_T)$
  & Eqs.~\eqref{eq:sar-cs} and~\eqref{eq:sar-cs-critical} \\
II & Integral link, $\rho$ bounded & $n\ll n_T\ll n^{\frac{2\beta+3}{2\beta+1}}$
  (SL/TL; no SLP in SD/TD)
  & $(nn_T)^{-1/2}$
  & SL/TL: see~\eqref{eq:sar-ii} \\
III & Integral link, Beta$(a,1)$
  (phase at $a=2,3$)
  & $1\!<\!a\!<\!2$: as~II;\;
    $2\!<\!a\!<\!3$ or $a\!>\!3$:
    \newline $n^{\frac{2(\beta+1)}{2\beta+a}}\ll n_T
              \ll n^{\frac{2\beta+3}{2\beta+1}}$
  & see~\eqref{eq:iii-mild}--\eqref{eq:iii-spatial};
    $n_T^{-1}$, $(nn_T)^{-1/2}$, or
    $\bigl(n_T^{\frac{2\beta+1}{2\beta+3}}n^{-1}\bigr)^{\!\frac{1}{a-1}}
    \times n_T^{-\frac{2(\beta+1)}{2\beta+3}}$
  & Cors.~\ref{cor:sar-iii} and~\ref{cor:sar-iii-critical};
    Table~\ref{tab:iii-phase} \\
IV & Derivative link, Beta$(a,1)$ ($\beta>1$)
  & $a>2+\frac{3}{2(\beta-1)}$;
  $n^{\frac{2\beta-1}{2\beta+a}}\ll n_T\ll n^{\frac{2\beta-1}{2\beta+1}}$
  (at equality: the subpolynomial layer in Cor.~\ref{cor:sar-iv};
  no SLP below the threshold)
  & above threshold:
    $\bigl(n_T^{\frac{2\beta+1}{2\beta-1}}n^{-1}\bigr)^{\!\frac{1}{a-1}}$
    \newline $\times n_T^{-\frac{2(\beta-1)}{2\beta-1}}$;
    equality: \eqref{eq:iv-critical}
  & see~\eqref{eq:sar-iv} \\
\midrule
Neg. & Derivative link, $\rho$ \emph{bounded} & \emph{Never}
  & $\min(n^{-\frac{2(\beta-1)}{2\beta+1}},n_T^{-\frac{2(\beta-1)}{2\beta-1}})$
  & $\Theta(1)$ \\
\bottomrule
\end{tabular}
\end{table}

\paragraph{Mechanisms underlying the SLP.}
Each gain above is a division of labor, not a pooling of samples.
Under bounded density conditions, integration averages noise and permits
a parametric source-only risk.
The target observes $g$ directly and contributes information about its local
structure, yielding the joint rate $(nn_T)^{-1/2}$ in Result~II.  A density
singularity adds a second division. For $a>2$, source estimation of the
cumulative coefficients of $g$ is limited by the sparse observations near
zero; the target estimates this finite-dimensional component at rate $n_T^{-1}$.
If the singularity is strong enough, the two
samples also split space: the source carries the bulk and the target the
neighborhood of the singularity (Result~III).  The derivative link reverses the
first mechanism.  Differentiation amplifies noise, and on a window of scale $t$
a source observation of $f$ informs $\dot f$ only after a factor $t^2$.  The
source is then a weak channel for the target estimand.  Only a much stronger
singularity restores synergy (Result~IV), and under a bounded ratio nothing
does (Theorem~\ref{thm:deriv}).

\paragraph{Necessary conditions for SLP.}
Three ingredients recur.  (I) \emph{Difficulty asymmetry}: the two problems
differ in intrinsic difficulty, in effective sample size (Results~I and~III) or
in smoothness class (Results~II and~IV).  (II) \emph{Complementary
information}: each fills a gap the other cannot, spatially near and away from
the singularity (Results~I and~III), or through source information about
the running integral and direct target observations of $g$ (Result~II).  (III) \emph{Sample-size window in which SAR diverges}: the
SLP holds on the middle wedges of Figure~\ref{fig:phase}, always bracketed by
SD and TD, with the SAR attaining its largest order at the
equalizing scale in the displayed settings.

\paragraph{Threshold ordering across links.}
The singularity needed for a polynomial SLP window is ordered by the link: none for the
integral link (phase cuts at $a=2,3$), $a>2+\tfrac{1}{2\beta}$ for the identity
link, and $a>2+\tfrac{3}{2(\beta-1)}$ for the derivative link, the strictest.
Condition~(I) alone is not enough: the bounded-ratio derivative link satisfies
it and still admits no gain (Theorem~\ref{thm:deriv}).

\section{Numerical and Real-Data Experiments}
\label{sec:experiments}

We use simulations and an analysis of IO-VNBD to examine how the density
ratio, functional link, and relative sample sizes affect the gains from
combining source and target data.  Simulations compare source-only,
target-only, and joint estimators under controlled bounded and Beta-source
designs.  The real-data analysis compares integral and derivative links
between two observed signals from the same trajectories.  A separate
semi-synthetic study examines unbounded covariate shift using constructed
designs and functions with empirically calibrated observation errors
(Appendix~\ref{app:semi-synthetic-unbounded}).

\begin{figure}[t]
  \centering
  \includegraphics[width=\linewidth]{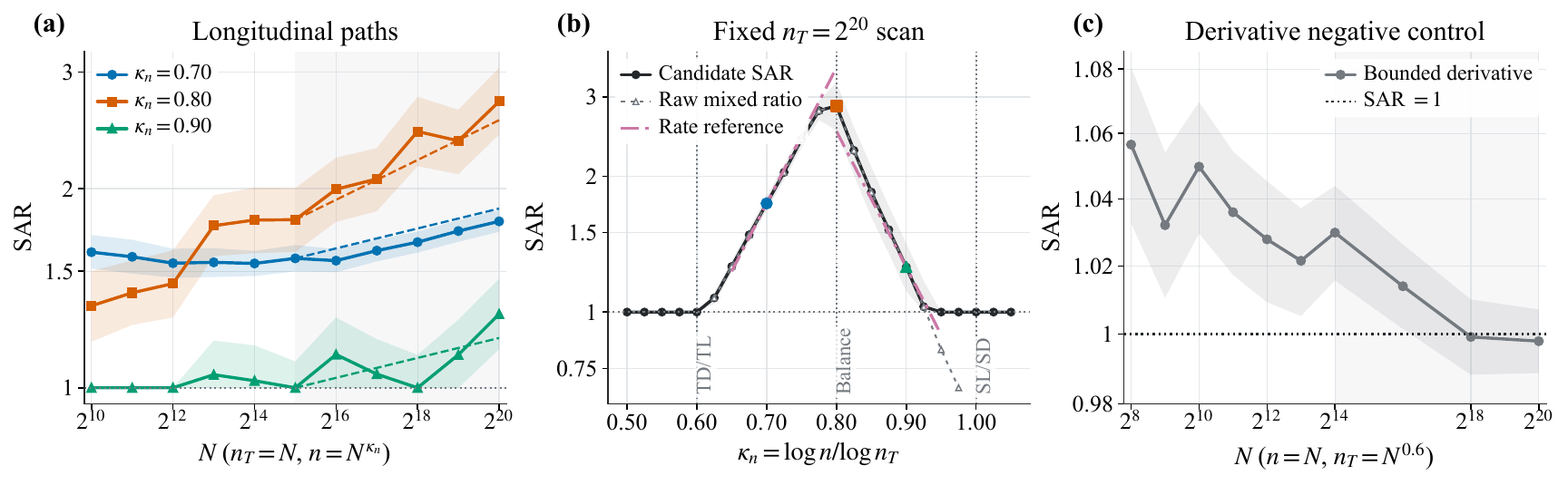}
  \caption{Integral--derivative comparison under uniform designs.
  (a) Integral candidate SAR for $n_T=N$, $n=N^{\kappa_n}$ with
  $\kappa_n=0.7,0.8,0.9$ and $\beta=1$.
  (b) Fixed-target scan at $n_T=2^{20}$; gray triangles show the raw mixed
  ratio only for $0.6<\kappa_n<1$.
  (c) Reversed derivative control for the same function pair, with $\beta=2$,
  $n=N$, $n_T=N^{0.6}$, and raw joint SAR.
  Bands are pointwise paired-bootstrap 95\% intervals. Dashed and dash--dot
  guides show theoretical slopes, vertical dotted lines mark theoretical
  boundaries and the equalizer, and shaded regions mark fitting windows.
  Full design and display details are in Appendix~\ref{app:sim:ii}.}
  \label{fig:sim:ii}
\end{figure}

\paragraph{Simulations under bounded designs.}
Figure~\ref{fig:sim:ii} compares the integral and derivative links under
uniform source and target designs, focusing on the effects of relative sample sizes
and the contrast between the two functional links.
Full experimental details are provided in Appendix~\ref{app:sim:ii}.

For the fixed truth in Figure~\ref{fig:sim:ii}, the balanced integral path
$\kappa_n=0.8$ has tail SAR slope $0.112$ $[0.065,0.159]$, versus the
theoretical reference $0.100$.
The fixed-target scan peaks at $\kappa_n=0.800$, while the derivative control
remains near one, ending at $0.998$ $[0.989,1.007]$.
Integral SAR uses the minimum of candidate mean risks, not an implemented
selector, so its lower bound of one is imposed by construction.
The scan's outer baseline follows prescribed switching and does not
independently validate the no-SLP regions or their boundaries.

\begin{figure}[!t]
  \centering
  \includegraphics[width=\linewidth]{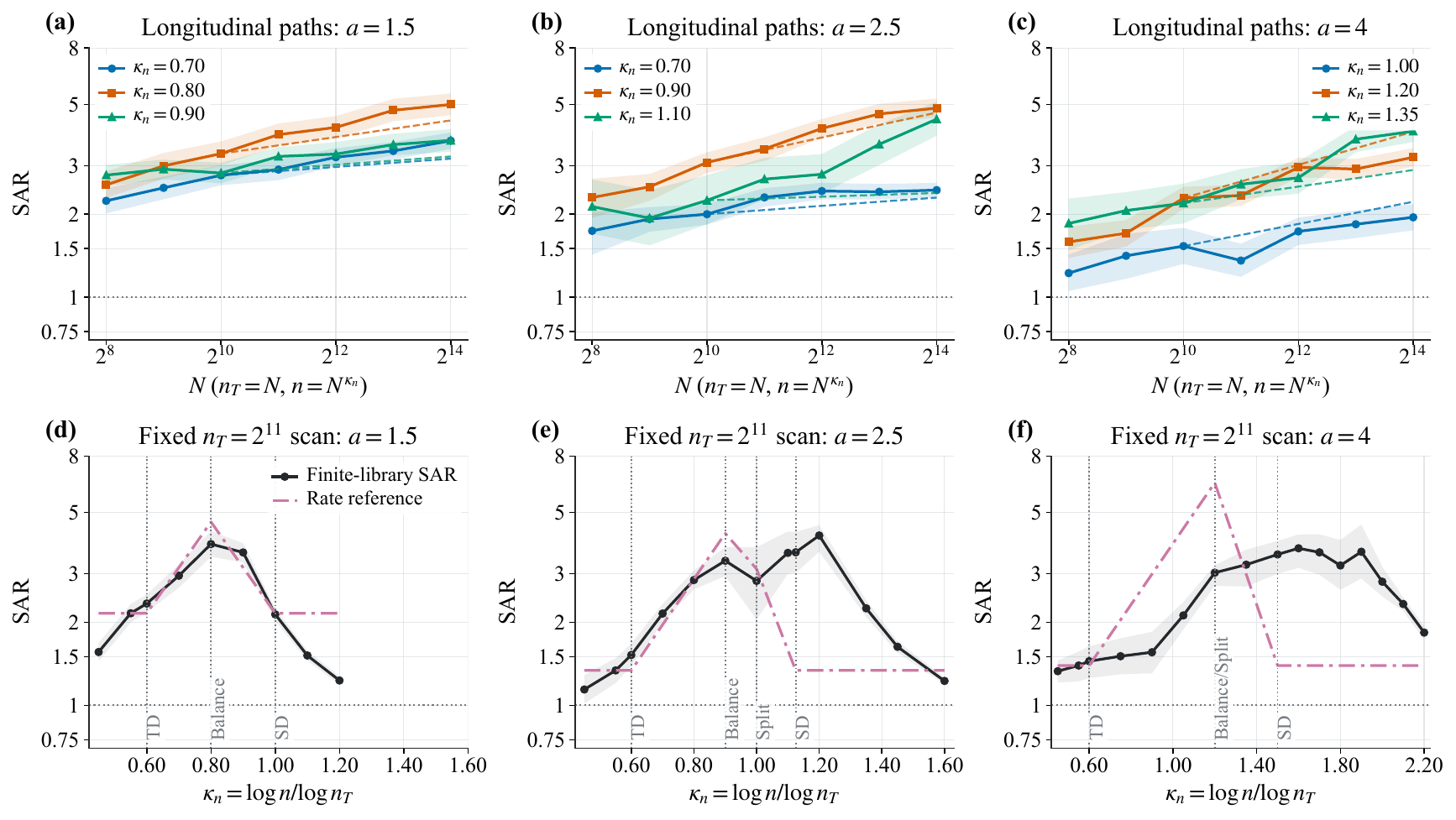}
  \caption{Regression experiments for the integral link under Beta-source designs
  with $a=1.5,2.5,4$, $\beta=1$, a fixed 64-dimensional representation,
  and cumulative bias bounds using the common smoothness radius $L=64$. Panels
  (a)--(c) use $n_T=N$ and $n=N^{\kappa_n}$ and show SAR for predeclared
  mechanism-specific truths; same-color dashed segments are theoretical
  tail-rate shapes aligned at the first fitted scale.  Panels (d)--(f) fix
  $n_T=2^{11}$ and show finite-library SAR against
  $\kappa_n=\log n/\log n_T$; purple dash--dot curves are aligned theoretical
  shapes and vertical dotted lines mark the TD boundary, equalizer,
  mechanism split, and SD boundary when distinct.  Panels (e,f) use paired
  trajectories on grids ending at $\kappa_n=1.60,2.20$, respectively.
  Bands are 95\% bootstrap intervals and horizontal dotted lines mark SAR$=1$.
  Full design and quantitative comparisons are in Appendix~\ref{app:sim:iii}.}
  \label{fig:sim:iii-noncritical}
\end{figure}

\paragraph{Simulations under Beta-source designs.}
We also examine the integral link with a uniform target and
Beta$(a,1)$ source designs at $a=1.5,2.5,4$, covering the three
parameter ranges separated by $a=2$ and $a=3$.
With $\beta=1$, we compare source-only, target-only, and joint
coefficient estimators in a fixed 64-dimensional representation,
using 400 independent regression repetitions per configuration.
Figure~\ref{fig:sim:iii-noncritical} shows SAR along joint sample-size
paths for prespecified truths and finite-library SAR in fixed-target scans.

The experiments show gains from combining the samples over parts of the
tested ranges, and the scans exhibit declining gains for source-rich
allocations. Several SAR slopes differ from their asymptotic references,
and the scan maxima at $a=2.5,4$ are displaced from the theoretical equalizers.
These fixed-truth and finite-library risks concern the specified estimators
rather than minimax optimization over the full H\"older class.
The experiments therefore illustrate finite-sample behavior without
establishing the full-class minimax rates or exact SLP boundaries.
Implementation details and quantitative comparisons are provided in
Appendix~\ref{app:sim:iii}. Expected-risk calculations at $a=2$ and $a=3$,
with selected regression Monte Carlo checks, are also reported there.

\paragraph{Real-data analysis.}
The public IO-VNBD benchmark \citep{onyekpe2021iovnbd} records vehicle and
smartphone measurements on public-road trajectories.  We use vehicle yaw rate
recorded through the controller area network (CAN) and smartphone GPS heading,
synchronized at 10~Hz.  Yaw rate measures
the vehicle's turning rate, while GPS heading records its direction of
travel, providing a physical basis for the functional links studied here.
After heading unwrapping, sign-and-scale calibration, and normalization,
the heading change $g$ is approximately the integral of the yaw-rate signal
$f$ on the audited windows, $g(x)\approx\int_0^x f(u)\,du$.
Taking yaw rate as the source response and heading change as the target
response gives the integral task; exchanging their roles gives the
derivative task.  The dataset therefore permits comparison of the two links
using the same pair of observed signals, without constructing synthetic
responses.  Data selection and calibration are detailed in
Appendix~\ref{app:real-data}.

Across 21 audited windows from 11 signal-distinct tracks, the terminal
track-geometric SAR is $1.385$ $[1.126,1.699]$ on the $\kappa=0.8$
integration path and $0.954$ $[0.930,0.976]$ for safeguarded differentiation
(Figure~\ref{fig:real-asymmetry}(b) in Appendix~\ref{app:real-data}).
These exploratory 95\% track-bootstrap
intervals condition on the selected windows and do not account for
window-selection uncertainty.
Finite-range rate discrepancies remain, and these illustrations do not
establish minimax rates.

\section{Discussion and Conclusion}
\label{sec:discussion}

Under a known functional link, synergy arises when the source and target
samples provide complementary information about the target regression
function, allowing joint estimation to improve on both individual minimax
rates for suitable relative sample sizes.
When both design densities are bounded above and away from zero, integration
admits SLP for suitable relative sample sizes, whereas differentiation does not.
Under Beta-source designs, the two links have different singularity requirements
and sample-size windows for SLP.  The minimax characterizations identify
the rate gains available from combining the samples, together with the
regimes in which the joint risk has the better individual rate.

\paragraph{Limitations.}
Only one-dimensional regression is studied. The rate-optimal estimators
require knowledge of the source density and $(a,\beta)$, so the adaptation
issue remains open, as do extensions to higher-order links and multivariate
singularities. Only asymptotic results on convergence rates are provided. 

\label{main-text-end}

\bibliographystyle{plainnat}
\bibliography{ref}

@inproceedings{pathak2022new,
  title={A new similarity measure for covariate shift with applications to nonparametric regression},
  author={Pathak, Reese and Ma, Cong and Wainwright, Martin},
  booktitle={International Conference on Machine Learning},
  pages={17517--17530},
  year={2022},
  organization={PMLR}
}

@article{cortes2010learning,
  title={Learning bounds for importance weighting},
  author={Cortes, Corinna and Mansour, Yishay and Mohri, Mehryar},
  journal={Advances in Neural Information Processing Systems},
  volume={23},
  year={2010}
}

@article{sugiyama2007direct,
  title={Direct importance estimation with model selection and its application to covariate shift adaptation},
  author={Sugiyama, Masashi and Nakajima, Shinichi and Kashima, Hisashi and Buenau, Paul and Kawanabe, Motoaki},
  journal={Advances in neural information processing systems},
  volume={20},
  year={2007}
}

@article{kpotufe2021marginal,
  title={Marginal singularity and the benefits of labels in covariate-shift},
  author={Kpotufe, Samory and Martinet, Guillaume},
  journal={The Annals of Statistics},
  volume={49},
  number={6},
  pages={3299--3323},
  year={2021},
  publisher={Institute of Mathematical Statistics}
}

@article{shimodaira2000improving,
  title={Improving predictive inference under covariate shift by weighting the log-likelihood function},
  author={Shimodaira, Hidetoshi},
  journal={Journal of statistical planning and inference},
  volume={90},
  number={2},
  pages={227--244},
  year={2000},
  publisher={Elsevier}
}

@article{stone1982optimal,
  title={Optimal global rates of convergence for nonparametric regression},
  author={Stone, Charles J},
  journal={The annals of statistics},
  pages={1040--1053},
  year={1982},
  publisher={JSTOR}
}

@article{reeve2021adaptive,
  title={Adaptive transfer learning},
  author={Reeve, Henry WJ and Cannings, Timothy I and Samworth, Richard J},
  journal={arXiv preprint arXiv:2106.04455},
  year={2021}
}

@incollection{torrey2010transfer,
  title={Transfer learning},
  author={Torrey, Lisa and Shavlik, Jude},
  booktitle={Handbook of research on machine learning applications and trends: algorithms, methods, and techniques},
  pages={242--264},
  year={2010},
  publisher={IGI global}
}

@article{mansour2009domain,
  title={Domain adaptation: Learning bounds and algorithms},
  author={Mansour, Yishay and Mohri, Mehryar and Rostamizadeh, Afshin},
  journal={arXiv preprint arXiv:0902.3430},
  year={2009}
}

@article{cai2021transfer,
  title={Transfer learning for nonparametric classification: Minimax rate and adaptive classifier},
  author={Cai, T Tony and Wei, Hongji},
  journal={The Annals of Statistics},
  volume={49},
  number={1},
  pages={100--128},
  year={2021},
  publisher={Institute of Mathematical Statistics}
}

@article{tian2021transfer,
  title={Transfer Learning under High-dimensional Generalized Linear Models},
  author={Tian, Ye and Feng, Yang},
  journal={arXiv preprint arXiv:2105.14328},
  year={2021}
}

@misc{zhou2022synergistic,
  author        = {Zhou, Ling and Yang, Yuhong},
  title         = {On a Synergistic Learning Phenomenon in Nonparametric
                   Domain Adaptation},
  year          = {2025},
  howpublished  = {arXiv:2511.17009},
  eprint        = {2511.17009},
  archivePrefix = {arXiv},
  primaryClass  = {stat.ML},
  url           = {https://arxiv.org/abs/2511.17009}
}

@article{pan2009survey,
  title={A Survey on Transfer Learning},
  author={Pan, Sinno Jialin and Yang, Qiang},
  journal={IEEE Transactions on Knowledge and Data Engineering},
  volume={22},
  number={10},
  pages={1345--1359},
  year={2009},
  publisher={IEEE}
}

@article{weiss2016survey,
  title={A Survey of Transfer Learning},
  author={Weiss, Karl and Khoshgoftaar, Taghi M and Wang, DingDing},
  journal={Journal of Big Data},
  volume={3},
  number={1},
  pages={9},
  year={2016},
  publisher={Springer}
}

@article{cavalier2008nonparametric,
  author    = {Cavalier, Laurent},
  title     = {Nonparametric statistical inverse problems},
  journal   = {Inverse Problems},
  volume    = {24},
  number    = {3},
  pages     = {034004},
  year      = {2008},
  doi       = {10.1088/0266-5611/24/3/034004}
}

@article{efromovich1994adaptive,
  author    = {Efromovich, Sam and Koltchinskii, Vladimir},
  title     = {On inverse problems with unknown operators},
  journal   = {IEEE Transactions on Information Theory},
  volume    = {47},
  number    = {7},
  pages     = {2876--2894},
  year      = {2001},
  doi       = {10.1109/18.959267}
}

@book{engl1996regularization,
  author    = {Engl, Heinz W and Hanke, Martin and Neubauer, Andreas},
  title     = {Regularization of Inverse Problems},
  publisher = {Kluwer Academic Publishers},
  year      = {1996}
}

@article{rice1984bandwidth,
  author    = {Rice, John A},
  title     = {Bandwidth choice for nonparametric regression},
  journal   = {The Annals of Statistics},
  volume    = {12},
  number    = {4},
  pages     = {1215--1230},
  year      = {1984}
}

@article{hall1990estimation,
  author    = {Hall, Peter and Marron, J S},
  title     = {On variance estimation in nonparametric regression},
  journal   = {Biometrika},
  volume    = {77},
  number    = {2},
  pages     = {415--419},
  year      = {1990}
}

@book{tsybakov2009introduction,
  author    = {Tsybakov, Alexandre B},
  title     = {Introduction to Nonparametric Estimation},
  publisher = {Springer},
  year      = {2009}
}

@inproceedings{tripuraneni2021provable,
  author    = {Tripuraneni, Nilesh and Jin, Chi and Jordan, Michael I},
  title     = {Provable meta-learning of linear representations},
  booktitle = {Proceedings of the 38th International Conference on Machine Learning},
  pages     = {10434--10443},
  year      = {2021}
}

@inproceedings{du2020fewshot,
  author    = {Du, Simon S and Hu, Wei and Kakade, Sham M and Lee, Jason D and Lei, Qi},
  title     = {Few-shot learning via learning the representation, provably},
  booktitle = {International Conference on Learning Representations},
  year      = {2021}
}

@article{li2022transfer,
  author    = {Li, Sai and Cai, T Tony and Li, Hongzhe},
  title     = {Transfer learning for high-dimensional linear regression: Prediction, estimation and minimax optimality},
  journal   = {Journal of the Royal Statistical Society Series B},
  volume    = {84},
  number    = {1},
  pages     = {149--173},
  year      = {2022}
}

@article{donoho1990minimax,
  author    = {Donoho, David L and Liu, Richard C},
  title     = {Geometrizing rates of convergence, {III}},
  journal   = {The Annals of Statistics},
  volume    = {19},
  number    = {2},
  pages     = {668--701},
  year      = {1991}
}

@article{nussbaum1999equivalence,
  author    = {Nussbaum, Michael},
  title     = {Minimax risk: Pinsker bound},
  journal   = {Encyclopedia of Statistical Sciences},
  year      = {1999}
}

@article{cai2024transfer,
  title     = {Transfer Learning for Nonparametric Regression: Non-asymptotic Minimax Analysis and Adaptive Procedure},
  author    = {Cai, T. Tony and Pu, Hongming},
  journal   = {arXiv preprint arXiv:2401.12272},
  year      = {2024}
}

@inproceedings{lin2024smoothness,
  title     = {Smoothness Adaptive Hypothesis Transfer Learning},
  author    = {Lin, Haotian and Reimherr, Matthew},
  booktitle = {Proceedings of the 41st International Conference on Machine Learning},
  series    = {Proceedings of Machine Learning Research},
  volume    = {235},
  pages     = {30286--30316},
  year      = {2024},
  url       = {https://proceedings.mlr.press/v235/lin24q.html}
}

@article{onyekpe2021iovnbd,
  author  = {Onyekpe, Uche and Palade, Vasile and Kanarachos, Stratis and
             Szkolnik, Alicja},
  title   = {{IO-VNBD}: Inertial and Odometry Benchmark Dataset for Ground
             Vehicle Positioning},
  journal = {Data in Brief},
  volume  = {35},
  pages   = {106885},
  year    = {2021},
  doi     = {10.1016/j.dib.2021.106885}
}

\clearpage
\appendix

\section{Function Class and Notation}
\label{app:notation}

With $r_\beta=\lceil\beta\rceil-1$ and $\alpha_\beta=\beta-r_\beta\in(0,1]$, the
H\"older class $\mcH(\beta,\kappa)$ collects the $f$ with
$\max_{j\leq r_\beta}\|f^{(j)}\|_\infty\leq\kappa$ whose $r_\beta$-th derivative
is $\alpha_\beta$-H\"older with constant $\kappa$; integer
smoothness $m$ thus means $C^{m-1,1}$.  Throughout, $b$ is a
bandwidth, the Beta shape $a$ is fixed, and $c_j,C_j$ are constants
independent of $(n,n_T)$.
For positive sequences, $A\preceq B$ means $A=O(B)$,
$A\succeq B$ means $B=O(A)$, and $A\asymp B$ means both.
We write $A\ll B$ for $A=o(B)$, i.e.\ $A/B\to0$, and $A\gg B$
for $B=o(A)$, along the sample-size sequence under consideration.
We use $\log_+(u):=\max\{\log u,0\}$ for $u>0$.

\section{The Identity Link: Covariate Shift}
\label{sec:covariate}
\label{sec:slp1}
\label{app:identity}

With $\mathcal{A}=\mathrm{Id}$ in~\eqref{eq:known-link} both domains regress the
same function, so the design pair is the only structural difference between
the domains, with density ratio $\rho=h_T/h$.  This identity-link setting is settled by
\citet{zhou2022synergistic}; we record it in SAR language as the reference
point for the non-identity links, and add the boundary case its analysis left
open.  Take a Beta$(a,1)$ source, $h(x)=ax^{a-1}$, a uniform target, and
$f\in\mcH(\beta,\kappa)$, so $\rho(x)\to\infty$ as $x\to0$ once $a>1$.
For $a>2+\tfrac{1}{2\beta}$, the SLP occurs precisely when
$n^{(2\beta+1)/(2\beta+a)}\ll n_T\ll n$;
there the minimax rate is $(n_T/n)^{1/(a-1)}n_T^{-2\beta/(2\beta+1)}$, faster
than either domain alone, the plane splits into the four regions SD/SL/TL/TD,
and the SAR attains its largest order at the equalizing scale
$n_T\asymp n^{(2\beta+1)^2/(2\beta(2\beta+a))}$.  At
$a=2+\tfrac{1}{2\beta}$, which is new here, the rate acquires a logarithmic
correction and the SLP survives only on the narrow band $n_T/n\to0$ with
$\log(n/n_T)=o(\log n)$.  Appendix~\ref{app:deferred:i} gives the full
theorem, the SAR formulas, and the sample-size regions.

\section{Additional Rates, SAR, and Phase Diagrams}
\label{app:deferred}

This appendix supplements the main-text summary in Table~\ref{tab:summary}
with additional rate and SAR formulas, sample-size-region descriptions
for Results~I and IV, the SAR formulas for Result~III at $a=2$ and $a=3$, and its
phase table.
The density conditions, mechanisms, minimax rates, and SAR formulas
for Result~III when $1<a<2$, $2<a<3$, or $a>3$ are presented in Section~\ref{sec:slp3}.
The SAR formulas for Result~IV are given in Corollary~\ref{cor:sar-iv}.
The displayed rate and SAR formulas use the same region labels as the main
text: SL and TL denote sequences on which the SAR diverges, with strict
outer endpoints.  The intermediate risk ranges at the identity- and
derivative-link thresholds are stated separately from their SLP conditions.

For the Result~III comparisons in this appendix,
write $n_T^{\mathrm{TD}}\asymp n^{(2\beta+3)/(2\beta+1)}$ and, for $a>2$,
$n_T^{\mathrm{cum}}\asymp n^{2(\beta+1)/(2\beta+a)}$.
For $a>3$, define the spatial-transition scale and rate benchmark by
\begin{equation}\label{eq:iii-rspat}
    n_T^{\mathrm{sp}}:=n^{\frac{2\beta+3}{2\beta+a}},
    \qquad
    R_{\mathrm{spat}}(n,n_T;a,\beta)
    :=\left(n_T^{\frac{2\beta+1}{2\beta+3}}n^{-1}\right)^{\frac{1}{a-1}}
    n_T^{-\frac{2(\beta+1)}{2\beta+3}}.
\end{equation}
$R_{\mathrm{spat}}$ is a rate benchmark rather than an exact
risk. It equals $n_T^{-1}$ at $n_T=n_T^{\mathrm{sp}}$ and matches the
target-only rate at $n_T\asymp n_T^{\mathrm{TD}}$.

\subsection{Full statement, SAR, and sample-size regions for Result~I}
\label{app:deferred:i}

The target-only rate is $n_T^{-2\beta/(2\beta+1)}$; the source-only rate is
$n^{-2\beta/(2\beta+1)}$ for $a<2+\tfrac{1}{2\beta}$,
carries a $\log(en)$ factor at equality, and slows to
$n^{-(2\beta+1)/(2\beta+a)}$ above it.

\begin{theorem}[Minimax rate under covariate shift]
\label{thm:cs}
Let $h(x)=ax^{a-1}$, $h_T\equiv1$, and
$f\in\mcH(\beta,\kappa)$.  For any $a>0$, $\beta>0$:
\begin{itemize}[leftmargin=*,topsep=2pt,itemsep=2pt,parsep=0pt]
\item If $a< 2+\tfrac{1}{2\beta}$: $\;\Rht(\mcH(\beta,\kappa);\,n,n_T)
      \asymp(n+n_T)^{-2\beta/(2\beta+1)}$.\; No SLP.
\item If $a=2+\tfrac{1}{2\beta}$: the rate is
      $\mathfrak C_\beta(n,n_T)$ of~\eqref{eq:cs-critical}, stated in
      Appendix~\ref{app:deferred:i-critical}.
\item If $a>2+\tfrac{1}{2\beta}$:
      \begin{equation}\label{eq:cs}
        \Rht(\mcH(\beta,\kappa);\,n,n_T)\;\asymp\;
        \begin{cases}
          n^{-\frac{2\beta+1}{2\beta+a}},
          & n_T\preceq n^{\frac{2\beta+1}{2\beta+a}}
            \quad\text{(SD)},\\[4pt]
          \bigl(n_T n^{-1}\bigr)^{\frac{1}{a-1}}n_T^{-\frac{2\beta}{2\beta+1}},
          & n^{\frac{2\beta+1}{2\beta+a}}\ll n_T
            \preceq n^{\frac{(2\beta+1)^2}{2\beta(2\beta+a)}}
            \quad\text{(SL)},\\[4pt]
          \bigl(n_T n^{-1}\bigr)^{\frac{1}{a-1}}n_T^{-\frac{2\beta}{2\beta+1}},
          & n^{\frac{(2\beta+1)^2}{2\beta(2\beta+a)}}\preceq n_T
            \ll n \quad\text{(TL)},\\[4pt]
          n_T^{-\frac{2\beta}{2\beta+1}},
          & n_T\succeq n \quad\text{(TD)}.
        \end{cases}
      \end{equation}
\end{itemize}
The rates for $a\ne2+1/(2\beta)$ are from \citet{zhou2022synergistic}; at
$a=2+1/(2\beta)$, we sharpen the bound by specifying its logarithmic factor
(see Appendix~\ref{app:deferred:i-critical}).
\end{theorem}

The upper bounds in all stated parameter cases remain valid
for independent observations with conditionally mean-zero sub-Gaussian
noises and uniformly bounded conditional sub-Gaussian parameters.

\begin{corollary}[SAR under covariate shift]
\label{cor:sar-cs}
For $a>2+\tfrac{1}{2\beta}$, the SLP occurs if and only if
$n^{(2\beta+1)/(2\beta+a)}\ll n_T\ll n$, with the SAR given
in~\eqref{eq:sar-cs} below.
At $a=2+1/(2\beta)$ the SLP occurs exactly when $n_T/n\to0$ with
$\log(n/n_T)=o(\log n)$.  For instance, $n_T=n/(\log n)^k$ qualifies while
$n_T=n^\tau$, $\tfrac{2\beta}{2\beta+1}<\tau<1$, does not.
The explicit SAR at equality is given in
Appendix~\ref{app:deferred:i-critical}.
\end{corollary}

In that window,
\begin{equation}\label{eq:sar-cs}
  \SAR\asymp
  \begin{cases}
    \bigl(n_T\,n^{-\frac{2\beta+1}{2\beta+a}}\bigr)^{\frac{2\beta}{2\beta+1}-\frac{1}{a-1}},
    & n^{\frac{2\beta+1}{2\beta+a}}\ll n_T
      \preceq n^{\frac{(2\beta+1)^2}{2\beta(2\beta+a)}}
      \quad\text{(SL)},\\[4pt]
    (n n_T^{-1})^{\frac{1}{a-1}},
    & n^{\frac{(2\beta+1)^2}{2\beta(2\beta+a)}}\preceq n_T
      \ll n \quad\text{(TL)}.
  \end{cases}
\end{equation}

For $a>2+1/(2\beta)$, three boundaries split the $(n,n_T)$-plane into four
sample-size regions (SSR), in which source and target play different roles
\citep{zhou2022synergistic}:
\begin{enumerate}[leftmargin=*, topsep=2pt, itemsep=1pt]
\item \textbf{Source Dominating (SD)}
  ($n_T\preceq n^{(2\beta+1)/(2\beta+a)}$):
  $n_T$ is too small to help near the singularity $x=0$;
  the rate equals the source-only rate and $\SAR=\Theta(1)$.
\item \textbf{Source Leading (SL)}
  ($n^{(2\beta+1)/(2\beta+a)}\mathrel{\ll} n_T
  \preceq n^{(2\beta+1)^2/(2\beta(2\beta+a))}$):
  source alone is still better than target alone, but target data
  assist near $x=0$ where $h$ vanishes; $\SAR\to\infty$ and
  its rate increases with $n_T$.
\item \textbf{Target Leading (TL)}
  ($n^{(2\beta+1)^2/(2\beta(2\beta+a))}\preceq n_T\mathrel{\ll} n$):
  target alone is already better; the larger source sample assists
  away from $x=0$ where $h$ is not small; $\SAR\asymp(n/n_T)^{1/(a-1)}$
  has a decreasing rate that returns to order one at the TD boundary.
\item \textbf{Target Dominating (TD)} ($n_T\succeq n$):
  source data no longer improve the rate; $\SAR=\Theta(1)$.
\end{enumerate}
SLP holds only in SL and TL.  The equalizing size
$n_T\asymp n^{(2\beta+1)^2/(2\beta(2\beta+a))}$ is where the two
individual rates match and the SAR attains its largest order.

\subsection{SAR for the identity link at \texorpdfstring{$a=2+1/(2\beta)$}{a=2+1/(2 beta)}}
\label{app:deferred:i-critical}

At $a=2+1/(2\beta)$ the common-function rate
$\mathfrak C_\beta(n,n_T)$, which gives the \emph{order} of the joint risk when
both samples observe the same function, is defined for $n,n_T\geq1$ by
      \begin{equation}\label{eq:cs-critical}
       \begin{aligned}
       \mathfrak C_\beta(n,n_T)&:=
       \begin{cases}
        n^{-\frac{2\beta}{2\beta+1}}\log(en),
          &n_T\leq n^{\frac{2\beta}{2\beta+1}},\\[4pt]
        n^{-\frac{2\beta}{2\beta+1}}\{1+\log(n n_T^{-1})\},
          &n^{\frac{2\beta}{2\beta+1}}<n_T<n,\\[4pt]
        n_T^{-\frac{2\beta}{2\beta+1}},&n_T\geq n.
        \end{cases}\\[5pt]
       \Rht(\mcH(\beta,\kappa);\,n,n_T)
       &\asymp\mathfrak C_\beta(n,n_T).
       \end{aligned}
      \end{equation}

In the intermediate risk regime $n^{2\beta/(2\beta+1)}<n_T<n$ the SAR is
\begin{equation}\label{eq:sar-cs-critical}
 \SAR\asymp
 \begin{cases}
 \dfrac{\log(en)}{1+\log(n n_T^{-1})},
   &n^{\frac{2\beta}{2\beta+1}}<n_T
     \leq n\{\log(en)\}^{-\frac{2\beta+1}{2\beta}},\\[9pt]
 \dfrac{(n n_T^{-1})^{\frac{2\beta}{2\beta+1}}}{1+\log(n n_T^{-1})},
   &n\{\log(en)\}^{-\frac{2\beta+1}{2\beta}}<n_T<n.
 \end{cases}
\end{equation}
The crossover follows by equating the individual rates
$n^{-2\beta/(2\beta+1)}\log(en)$ and $n_T^{-2\beta/(2\beta+1)}$.  Outside this intermediate regime,
$\SAR\asymp1$.  The intermediate risk regime is not itself an SLP
window: at this value of $a$, SLP occurs exactly when
$\log(n/n_T)=o(\log n)$ and $n_T/n\to0$, or equivalently
$z=n/n_T\to\infty$ with $\log z=o(\log n)$.
For example, $n_T=n/(\log n)^k$ gives SLP for every fixed $k>0$,
whereas $n_T=n^\tau$ with fixed $\frac{2\beta}{2\beta + 1}<\tau<1$ does not.

\subsection{SAR for Result~III at \texorpdfstring{$a=2$ and $a=3$}{a=2 and a=3}}
\label{app:deferred:iii-sar}

\begin{corollary}[SAR for Result~III at $a=2$ and $a=3$]
\label{cor:sar-iii-critical}
Under the setting of Section~\ref{sec:iii-setting}, the SAR formulas
at $a=2$ and $a=3$ are as follows.
At $a=2$, the individual rates $(\log n)/n$ and
$n_T^{-2(\beta+1)/(2\beta+3)}$ cross at
\[
  n_T\asymp(n/\log n)^{(2\beta+3)/(2(\beta+1))}.
\]
In the SLP window,
\begin{equation}\label{eq:sar-iii-a2}
  \SAR\asymp
  \begin{cases}
    (\log n)\,n_T n^{-1},
      & n(\log n)^{-1}\ll n_T\preceq n \quad\text{(SL)},\\[6pt]
    (\log n)\,(n_T n^{-1})^{1/2},
      & n\preceq n_T\preceq
        \bigl(n(\log n)^{-1}\bigr)^{\frac{2\beta+3}{2(\beta+1)}}
        \quad\text{(SL)},\\[6pt]
    \bigl(n\,n_T^{-\frac{2\beta+1}{2\beta+3}}\bigr)^{1/2},
      & \bigl(n(\log n)^{-1}\bigr)^{\frac{2\beta+3}{2(\beta+1)}}
        \preceq n_T\ll n^{\frac{2\beta+3}{2\beta+1}}
        \quad\text{(TL)}.
  \end{cases}
\end{equation}
Thus the SLP window is $n/\log n\ll n_T\ll n^{(2\beta+3)/(2\beta+1)}$.

At $a=3$, the individual rates are
$n^{-2(\beta+1)/(2\beta+3)}$ and
$n_T^{-2(\beta+1)/(2\beta+3)}$.  In the SLP window,
\begin{equation}\label{eq:sar-iii-a3}
  \SAR\asymp
  \begin{cases}
    n_T n^{-\frac{2(\beta+1)}{2\beta+3}},
      & n^{\frac{2(\beta+1)}{2\beta+3}}\ll n_T
        \preceq n(\log n)^{-1} \quad\text{(SL)},\\[6pt]
    n^{\frac{1}{2\beta+3}}(\log n)^{-1},
      & n(\log n)^{-1}\preceq n_T\preceq n \quad\text{(SL)},\\[9pt]
    \dfrac{\bigl(n\,n_T^{-\frac{2\beta+1}{2\beta+3}}\bigr)^{1/2}}
      {1+\log_+\bigl(en\,n_T^{-\frac{2\beta+1}{2\beta+3}}\bigr)},
      & n\preceq n_T\ll n^{\frac{2\beta+3}{2\beta+1}}
        \quad\text{(TL)}.
  \end{cases}
\end{equation}
The SLP window is
$n^{2(\beta+1)/(2\beta+3)}\ll n_T\ll n^{(2\beta+3)/(2\beta+1)}$.
At $a=3$, logarithmic corrections alter the location of
the largest SAR order.  With $n_T=N$ held fixed, the individual rates match
at $n\asymp N$, whereas the SAR attains its largest order at
$n\asymp N\log N$.  For fixed $n$, Equation~\eqref{eq:sar-iii-a3} gives
a plateau in the SAR rate over $n/\log n\preceq n_T\preceq n$.
For both $a=2$ and $a=3$, $\mathrm{SAR}=\Theta(1)$
in the SD and TD regions.
The corresponding joint risk benchmarks are recorded in
Propositions~\ref{prop:iii-a2}--\ref{prop:iii-a3}, respectively.
\end{corollary}

\subsection{Phase-diagram table for Result~III}
\label{app:deferred:iii-table}

\begin{table}[H]
\centering
\caption{Result~III: SLP windows and synergistic rates
(integral link, Beta$(a,1)$).}
\label{tab:iii-phase}
\renewcommand{\arraystretch}{1.25}
\begin{tabular}{lll}
\toprule
Range of $a$ & SLP window & Synergistic rate(s) \\
\midrule
$1<a<2$
  & $n\ll n_T\ll n_T^{\mathrm{TD}}$
  & $(nn_T)^{-1/2}$ \\
$a=2$
  & $n/\log n\ll n_T\ll n_T^{\mathrm{TD}}$
  & $n_T^{-1}$ (thin), then $(nn_T)^{-1/2}$ \\
$2<a<3$
  & $n_T^{\mathrm{cum}}\ll n_T\ll n_T^{\mathrm{TD}}$
  & $n_T^{-1}$, then $(nn_T)^{-1/2}$ \\
$a=3$
  & $n^{2(\beta+1)/(2\beta+3)}\ll n_T\ll n_T^{\mathrm{TD}}$
  & \shortstack[l]{$n_T^{-1}$, then $(\log n)n^{-1}$,\\
      then log-corrected $(nn_T)^{-1/2}$} \\
$a>3$
  & $n_T^{\mathrm{cum}}\ll n_T\ll n_T^{\mathrm{TD}}$
  & $n_T^{-1}$, then $R_{\mathrm{spat}}(n,n_T;a,\beta)$ \\
\bottomrule
\end{tabular}
\end{table}

\subsection{Source-only rate for the derivative link}
\label{app:deferred:iv-src}

\begin{remark}[Source-only rate under Beta$(a,1)$ density]
\label{rem:iv-src-rate}
Under bounded source density, the source-only derivative rate is
$n^{-2(\beta-1)/(2\beta+1)}$.  Under the Beta source, the
same power holds for $a<2+3/\{2(\beta-1)\}$, acquires a logarithmic factor
at equality, and slows to $n^{-(2\beta-1)/(2\beta+a)}$ above it.
In the strong-singularity regime the combined rate~\eqref{eq:iv} strictly beats
both single-task benchmarks in the synergistic window, so the SLP and SAR
conclusions of Corollary~\ref{cor:sar-iv} remain valid.
\end{remark}

\subsection{The case \texorpdfstring{$a=2+3/\{2(\beta-1)\}$}{a=2+3/[2(beta-1)]} of Theorem~\ref{thm:iv}}
\label{app:deferred:iv-critical2}
\label{app:deferred:iv-critical}

Under the setting of Theorem~\ref{thm:iv}, if
$a=2+\tfrac{3}{2(\beta-1)}$, then
\begin{equation}\label{eq:iv-critical}
  \Rht(\dot f;\,n,n_T)\asymp
  \begin{cases}
    n^{-\frac{2(\beta-1)}{2\beta+1}}\log n,
      &n_T\preceq n^{\frac{2(\beta-1)}{2\beta+1}},\\[4pt]
    n^{-\frac{2(\beta-1)}{2\beta+1}}
      \Bigl\{1+\log_+\bigl(en\,n_T^{-\frac{2\beta+1}{2\beta-1}}\bigr)\Bigr\},
      &n^{\frac{2(\beta-1)}{2\beta+1}}\preceq n_T
        \preceq n^{\frac{2\beta-1}{2\beta+1}},\\[4pt]
    n_T^{-\frac{2(\beta-1)}{2\beta-1}},
      &n_T\succeq n^{\frac{2\beta-1}{2\beta+1}}.
  \end{cases}
\end{equation}

At $n_T\asymp n^{(2\beta-1)/(2\beta+1)}$, the
logarithmic factor is of constant order and the middle rate matches
the target-only rate up to constants.
For the SAR, write the common logarithmic factor as
\[
  \Lambda_D=1+\log_+\bigl(en\,n_T^{-\frac{2\beta+1}{2\beta-1}}\bigr).
\]
In the intermediate sample-size range,
\[
\SAR\asymp
\begin{cases}
  \log n/ \Lambda_D,
    &n^{\frac{2(\beta-1)}{2\beta+1}}\preceq n_T\preceq
      n^{\frac{2\beta-1}{2\beta+1}}
      (\log n)^{-\frac{2\beta-1}{2(\beta-1)}},\\[9pt]
  n^{\frac{2(\beta-1)}{2\beta+1}}
    n_T^{-\frac{2(\beta-1)}{2\beta-1}}\big/ \Lambda_D,
    &n^{\frac{2\beta-1}{2\beta+1}}
      (\log n)^{-\frac{2\beta-1}{2(\beta-1)}}
      \preceq n_T\preceq n^{\frac{2\beta-1}{2\beta+1}}.
\end{cases}
\]
The two branches meet where the individual risks have the same order.
This intermediate risk range is not itself an SLP window: the SAR diverges
only when $1\ll n^{(2\beta-1)/(2\beta+1)}n_T^{-1}=n^{o(1)}$.
Thus any diverging SAR at equality is subpolynomial, and there is no
fixed polynomial-exponent SLP window.

\subsection{Why the integration--differentiation asymmetry holds}
\label{app:deferred:asym}

Integration \emph{raises} smoothness ($f\in\mcH(\beta)\Rightarrow g\in\mcH(\beta+1)$)
and, under bounded density conditions, permits a parametric source-only
rate despite $g$ being nonparametric.
Differentiation \emph{lowers} smoothness and amplifies source noise.
Under bounded density conditions, its joint minimax rate is the
minimum of the two individual rates in Theorem~\ref{thm:deriv}.
Thus target data can improve upon source-only derivative estimation,
but not yield a rate-level gain over the better individual experiment.

\section{Additional Simulation Results}
\label{app:simulation}

We illustrate selected transfer scenarios from Sections~\ref{sec:slp2}--\ref{sec:slp4}.  The
experiments examine rate mechanisms, sample-allocation effects, and the
limitations of finite-range comparisons with minimax theory.

\paragraph{Common experimental setting.}
Design points are sampled independently on the unit cube of the corresponding
dimension.  All regression experiments use independent Gaussian noise with
standard deviation $0.25$.  Every Monte Carlo configuration in this
simulation section and its diagnostics uses 400 independent
repetitions, including regression experiments, surrogate
loss simulations, and separate Monte Carlo diagnostics.  Analytic-risk
calculations have no Monte Carlo repetition count.  Unless a subsection
specifies otherwise, we approximate integrated risk on a 4096-point grid
and construct paired percentile-bootstrap 95\% intervals from 5000
resamples; bootstrap resamples are distinct from simulation repetitions.
Slopes use the largest four sample sizes unless another tail length is stated.

For the synthetic regression experiments, $\beta$ denotes the
source-function class index used for tuning and rate
references.  The target class index is $\beta+1$ for the integral link,
$\beta-1$ for the derivative link.
Individual test functions may be smoother than these class requirements.
The surrogate experiments specify their coefficient classes separately.

\paragraph{Risk interpretation.}
We distinguish raw-regression Monte Carlo, surrogate-model Monte Carlo,
and analytic-risk calculations.  Monte Carlo is abbreviated as MC.
Raw-regression Monte Carlo generates design points and responses
and applies the regression estimators.  Surrogate-model Monte Carlo
samples losses induced by a separately specified observation model and
estimator, possibly by exact aggregation of coefficient errors.
Analytic-risk calculations evaluate expected losses without Monte Carlo
sampling and may concern either a regression estimator or a surrogate.
Each experiment evaluates specified estimators at fixed truths or over
finite function libraries. The minimax formulas serve as asymptotic
benchmarks for these risks and their SARs.  Bootstrap bands quantify pointwise
Monte Carlo precision under the stated model and risk definition.
Discrepancy between the surrogate and regression models is outside their
scope. Slopes, normalized risks, and scan maxima describe behavior over
the displayed sample sizes, while conclusions about minimax optimality
and asymptotic SLP boundaries rely on the theoretical analysis.
A bounded SAR allows constant-factor gains while indicating no
asymptotic rate improvement.

\paragraph{Display conventions.}
Sample-size and SAR axes are logarithmic; normalized-risk and relative-risk
diagnostics use linear vertical axes.  We write $\kappa_n=\log n/\log n_T$
for target-anchored scans and $\kappa_T=\log n_T/\log n$ for source-anchored scans.
Markers identify Monte Carlo estimates.  Colored dotted curves give expected
risks or SAR formed from those risks for the specified estimators, without
sampling bands.  Rate references instead give theoretical order shapes.
Unless stated otherwise, same-color dashed curves give longitudinal
references and purple dash--dot curves give scan references.  Rate curves
are vertically aligned at the first fitted point because the theorems
determine rates only up to constants.  Light-gray regions mark slope windows.
Horizontal gray dotted lines mark SAR$=1$.  Vertical gray dotted lines mark
theoretical balances or boundaries,
orange dashed lines mark channel crossings computed from the displayed
risks, and gray dashed vertical lines mark estimator switches.
Experiment-specific reference normalizations are stated below.

\subsection{Bounded integral link and the bounded-derivative control}
\label{app:sim:ii}

\paragraph{Design.}
The integral-link experiments use $\beta=1$, so the target class has
smoothness index $\beta+1=2$.  Let
\[
 f(x)=0.7|x-0.4|,
 \qquad
 g(x)=0.35\{0.4^2+\operatorname{sgn}(x-0.4)|x-0.4|^2\},
\]
so that $g'=f$.  Both designs are uniform.  On the longitudinal paths,
$n_T=N$ and $n=\operatorname{round}(N^{\kappa_n})$, with
$\kappa_n\in\{0.7,0.8,0.9\}$ and $N=2^{10},\ldots,2^{20}$.  The source
estimator is the empirical IPW integral.  The target estimator is degree-one
local polynomial regression with bandwidth $1.1N^{-1/5}$.  The mixed
estimator applies the count-based target correction with
$b=1\wedge(n/n_T)^{1/2}$.

The transverse scan fixes $n_T=2^{20}$ and sets
$n=\operatorname{round}(n_T^{\kappa_n})$ on the grid
$\kappa_n=0.5,0.525,\ldots,1.05$.

\paragraph{Risk diagnostics and prescribed switching.}
Both integral-link panels report a fixed-truth candidate-risk diagnostic.
Let $\widehat R_E(g)$ denote estimator $E$'s mean loss over the Monte Carlo
replicates at the displayed truth, and define
\[
 \widehat R_{\mathrm{cand}}(g)
 =\min\{\widehat R_S(g),\widehat R_T(g),\widehat R_B(g)\},
 \qquad
 \widehat{\mathrm{SAR}}_{\mathrm{cand}}
 =\frac{\widehat R_S(g)\wedge\widehat R_T(g)}
 {\widehat R_{\mathrm{cand}}(g)}.
\]
The candidate minimum is taken \emph{after} averaging each estimator's
losses, not separately within each Monte Carlo replicate. This is an oracle
comparison of the candidate mean risks.
The minimum and numerator are recomputed inside each paired bootstrap draw.
Nevertheless, $\widehat{\mathrm{SAR}}_{\mathrm{cand}}
=\max\{1,(\widehat R_S\wedge\widehat R_T)/\widehat R_B\}$, so its
unit lower bound is imposed by construction.

On longitudinal paths, $B$ is the raw mixed estimator.  In the fixed-target
scan, the third candidate follows the predeclared theoretical branch rule
\[
 \widehat g_B=
 \begin{cases}
   \widehat g_T,&\kappa_n\leq0.6,\\
   \widehat g_M,&0.6<\kappa_n<1,\\
   \widehat g_S,&\kappa_n\geq1.
 \end{cases}
\]
The raw mixed estimator is evaluated on $0.6<\kappa_n<1$, where its risk
is shown by gray markers. Outside this interval, the prescribed
single-domain choice sets the candidate SAR to one.

\paragraph{Negative control and slope fitting.}
For the bounded-derivative negative control we reverse the operator
direction using the same function pair.  Source responses now observe
the primitive above, and target responses observe its derivative.
Thus the derivative-link setting has source index $\beta=2$ and target
index $\beta-1=1$, unlike the integral-link panels.
We take $n=N$ and
$n_T=\operatorname{round}(N^{0.6})$ for
$N=2^8,\ldots,2^{14},2^{16},2^{18},2^{20}$, and fit the operator-aware
joint local polynomial.  This path equalizes the two individual minimax
rates for this derivative setting.  The negative-control SAR uses the
actual joint-estimator risk, without a candidate-minimum floor.

The three integral longitudinal fits use the largest \emph{six} scales,
$N=2^{15},\ldots,2^{20}$.  The negative control uses four,
$N=2^{14},2^{16},2^{18},2^{20}$.  Fixed-target slopes use the six grid
points in each of $[0.65,0.775]$ and $[0.825,0.95]$ for $\kappa_n$,
regressing log SAR on $\log n$; the reference slopes $\pm1/2$ are not
slopes with respect to $\kappa_n$ itself.

\paragraph{Results.}
Figure~\ref{fig:sim:ii}(a) compares the three longitudinal paths.
Their SAR-slope intervals in Table~\ref{tab:sim:summary} each contain
the corresponding theoretical exponent.  The largest point estimate occurs
on the balanced path, consistent with the theoretical equalizer $\kappa_n=0.8$.
At $\kappa_n=0.9$, two of the six tail points hit the candidate SAR floor;
the raw mixed ratio instead has slope $0.057$ $[0.009,0.109]$ on the same
points.  The candidate absolute-risk slopes for $\kappa_n=0.7,0.8$
lie somewhat above their theoretical decay exponents.

Figure~\ref{fig:sim:ii}(b) traverses the phase diagram at fixed $n_T$.
Within the actually evaluated mixed region, candidate SAR rises and then
falls, with its grid maximum at $\widehat\kappa_n=0.800$ and a gridwise
bootstrap interval $[0.775,0.800]$.  The estimated log--log slopes are $0.462$
$[0.442,0.484]$ in the interior of the TL region and $-0.501$
$[-0.567,-0.421]$ in the SL region, compared with $1/2$ and $-1/2$.
The left interval excludes $1/2$; the right estimate includes the
candidate floor near $\kappa_n=0.95$.  On the same right-hand points,
the raw mixed ratio has slope $-0.580$ $[-0.671,-0.487]$.
The peak SAR is $2.86$ $[2.70,3.23]$, with the maximum recomputed over the
grid in each bootstrap draw, rather than a pointwise interval at the
estimated peak.  The theoretical value
$n_T^{1/10}=4$ is an order benchmark because the rate statement suppresses
multiplicative constants.  Near $\kappa_n=1$, the raw mixed estimator is
less competitive than source-only in finite samples. The candidate minimum
therefore returns to one before the prescribed SD boundary.

The bounded-derivative control remains near the no-gain baseline in
Figure~\ref{fig:sim:ii}(c).  Its absolute-risk exponent in
Table~\ref{tab:sim:summary} is close to the theoretical value $0.4$.
The terminal SAR is $0.998$ $[0.989,1.007]$. The small negative finite-range
slope reflects the unsafeguarded joint estimator approaching the individual-risk
benchmark over the sampled range. These patterns are consistent with the
integration--differentiation asymmetry under bounded designs.

\paragraph{Reading Figure~\ref{fig:sim:ii}.}
The integral-link panels (a,b) use $\beta=1$; the reversed derivative
  control (c) uses $\beta=2$.
  (a) Candidate SAR along
  $n_T=N$ and $n=N^{\kappa_n}$ for $\kappa_n=0.7,0.8,0.9$.  Dashed segments have
  theoretical slopes $0.05,0.10,0.05$, aligned over the six-point tail
  $N=2^{15},\ldots,2^{20}$.  (b) The fixed-$n_T$ cross-section.
  The dark curve uses the minimum of candidate mean risks; gray triangles
  show the raw mixed ratio only where mixed was evaluated, $0.6<\kappa_n<1$.
  The outer unit baseline follows from prescribed single-domain switching.
  Purple dash--dot segments have log-SAR versus log-$n$ slopes $1/2$ and
  $-1/2$. Vertical dotted lines mark the theoretical boundaries and equalizer.
  (c) The bounded-derivative negative
  control on the equal-rate path $n=N$, $n_T=N^{0.6}$, extending
  to $N=2^{20}$, using raw joint SAR without a candidate floor.
  Bands are paired-bootstrap 95\% intervals.  Light-gray regions mark the
  fitting windows, using $2^{14},2^{16},2^{18},2^{20}$ in (c);
  horizontal dotted lines mark SAR$=1$.

\subsection{Unbounded integral link: regression paths and the cases \texorpdfstring{$a=2,3$}{a=2,3}}
\label{app:sim:iii}

\paragraph{Regression methods at $a=1.5,2.5,4$.}
The experiments at $a=1.5,2.5,4$ use fresh regression samples with
$\beta=1$ and a fixed coefficient family: two anchored cumulative polynomials
and 62 compactly supported zero-tail modes on five complete dyadic levels.
Reconstruction uses the population Gram matrix, with minimum eigenvalue
$0.125$ and condition number $15.15$; this matrix is never estimated from data.
Source coefficients use truncated Taylor/IPW and integration-by-parts
statistics, target coefficients use direct projection, and deterministic
inverse-upper-variance weights combine the channels.  The five predeclared
truths---linear cumulative, boundary cumulative, interior tree, boundary tree,
and mixed tree---satisfy $\|g\|_\infty,\|g'\|_\infty\leq0.8$,
$\operatorname{Lip}(g')\leq L=64$, conditional second moments at most one,
and analysis-coefficient envelopes with constant one.  Cumulative cutoffs
minimize a variance-plus-bias bound with
$B_\nu(\epsilon)=L\sup_{[0,\epsilon]}|K_\nu|\epsilon^2$; the joint criterion
includes the squared source weight on this bias term.  Independent piecewise
quadrature checks every truth's actual truncation bias at every candidate
cutoff.  The exact finite-span loss is
$(\widehat\theta-\theta)^\top G^{-1}(\widehat\theta-\theta)$,
where $G$ is the population analysis Gram matrix.

\paragraph{Design at $a=1.5,2.5,4$.}
To match the notation of Section~\ref{sec:slp2}, every longitudinal path at $a=1.5,2.5,4$ uses
the target-anchored convention
\[
 n_T=N,\qquad n=N^{\kappa_n},\qquad N=2^8,2^9,\ldots,2^{14}.
\]
For $a=1.5$ we take $\kappa_n=0.70,0.80,0.90$; for $a=2.5$ we take
$\kappa_n=0.70,0.90,1.10$; and for $a=4$ we take
$\kappa_n=1.00,1.20,1.35$.  At $a=4$, the three paths lie in the
spatial branch, at the spatial--level split and equalizer, and in the level
branch, respectively.  At $a=2.5$, the level path lies inside the narrow
interval $1<\kappa_n<1.125$, with theoretical SAR exponent $1/45$.
The path-specific stress truth is fixed before simulation: the interior
tree for integral paths, the cumulative-boundary truth for level paths,
the boundary tree for the spatial path, and the mixed tree at the $a=4$
equalizer. Each longitudinal curve uses a mechanism-specific truth, so
differences between curves reflect both the truth and the allocation.

The lower row instead fixes $n_T=2^{11}$ and scans
$\kappa_n=\log n/\log n_T$ through every theoretical breakpoint.
The $a=2.5$ and $a=4$ scans use grids ending at
$\kappa_n=1.60$ and $2.20$, respectively, to examine the source-rich
return.  Within each $a=2.5,4$ trajectory, source samples are nested across
allocations and the same target sample
is reused.  This preserves each allocation's marginal regression model
while coupling the scan points.  Here we report
$R_E^{\mathcal L}=\max_{g\in\mathcal L}
\mathbb E\|\widehat g_E-g\|_2^2$ and
$\min(R_S^{\mathcal L},R_T^{\mathcal L})/R_C^{\mathcal L}$; both the
library maximum and channel minimum are recomputed in every bootstrap draw.
For these coupled scans, cross-allocation comparisons resample complete
trajectories to preserve their dependence.
Longitudinal slopes use the largest five values of~$N$.

\paragraph{Results at $a=1.5,2.5,4$.}
Figure~\ref{fig:sim:iii-noncritical} in the main text presents these regression experiments.
Panels (a)--(c) distinguish finite-sample
synergy from agreement with the theoretical SAR exponents.
At $a=1.5$, the empirical slopes are $0.106$ $[0.070,0.141]$,
$0.148$ $[0.102,0.192]$, and $0.093$ $[0.051,0.135]$, versus
references $0.050,0.100,0.050$.  The equalizer path has the largest point
estimate, but the reference exponents lie below the corresponding
intervals, so the agreement is qualitative.  At $a=2.5$, the integral-rise
and equalizer slopes are $0.065$ $[0.032,0.098]$ and $0.173$
$[0.135,0.211]$, compatible with $0.050$ and $0.150$.
The level path has SAR slope $0.238$ $[0.157,0.317]$, well above
$1/45\simeq0.022$.  Its mixed-risk slope, $-1.033$ $[-1.111,-0.954]$,
is compatible with $n_T^{-1}$, but target-only risk remains smaller than
source-only risk at every sampled scale, unlike the asymptotic ordering on
this path. The empirical SAR therefore compares the joint risk with the
target-only risk throughout this range.

At $a=4$, the spatial, equalizer, and level paths have SAR slopes
$0.113$ $[0.053,0.173]$, $0.130$ $[0.092,0.167]$, and $0.227$
$[0.164,0.289]$, versus $0.133,0.200,0.100$.
Only the spatial path is compatible with its SAR reference; the equalizer
grows more slowly and the level path more quickly than predicted.
Their normalized mixed-risk slopes, $0.009,0.032,-0.022$, all have
intervals containing zero, but target-only risk is still smaller on the
level path at every sampled scale.

The fixed-$n_T$ scans in panels (d)--(f) exhibit a source-rich
decline, while retaining finite-sample fluctuations and displaced maxima.
For $a=1.5$, the largest sampled SAR occurs
at the theoretical equalizer $\kappa_n=0.8$, with value $3.84$
$[3.46,4.28]$.  For $a=2.5$, the largest sampled value is $4.13$
$[3.61,4.50]$ at $\kappa_n=1.20$, falling to $1.23$ $[1.19,1.26]$
at $1.60$.  For $a=4$, it is $3.71$ $[3.15,4.16]$ at $1.60$,
falling to $1.84$ $[1.73,1.94]$ at $2.20$.
The predeclared right-side comparisons give SAR ratios
$\mathrm{SAR}(1.60)/\mathrm{SAR}(1.20)=0.297$ $[0.271,0.341]$
for $a=2.5$ and
$\mathrm{SAR}(2.20)/\mathrm{SAR}(1.90)=0.509$ $[0.401,0.644]$
for $a=4$, using paired trajectory bootstrap intervals.
At the two right endpoints, source-only risk is smaller than target-only
risk, with source/target ratios $0.085$ and $0.273$, respectively.
The sampled maxima remain to the right of the theoretical equalizers
$0.90$ and $1.20$, and neither
curve is monotone throughout each theoretical branch.

\paragraph{Regression design at $a=2,3$.}
For $a=2,3$, we use the random-design regression model with
$\beta=1$ and $g(0)=0$.
The fixed library consists of $f(x)=0$, $f(x)=0.3$, and $f(x)=0.2+0.4x$,
with $g(x)=\int_0^x f(t)\,dt$.  All satisfy
$\operatorname{Lip}(f)\leq H=64$ and
$\|f\|_\infty,\|g\|_\infty\leq0.8$.
A fixed infinite orthonormal piecewise-linear multiresolution system
replaces the fixed-dimensional analysis family used at $a=1.5,2.5,4$.
The truth libraries and estimators differ between the two constructions,
so comparisons with theory are made within each construction.
The implemented coefficient methods use compensated IPW or integration by
parts for source data, direct projection for target data, and deterministic
variance and bias bounds for combination and shrinkage.  The numerical risk
checks are given in Appendix~\ref{app:sim:iii-normalization}.  No Gaussian
prior is used.
We evaluate $R_E^{\mathcal L}=\max_{g\in\mathcal L}
\mathbb E\|\widehat g_E-g\|_2^2$ for each specified estimator and form
$\min(R_S^{\mathcal L},R_T^{\mathcal L})/R_C^{\mathcal L}$.

The longitudinal allocations are $n_T=N$ and $n=\operatorname{round}(N^{\kappa_n})$
for $\kappa_n=0.70,0.90,1.00$ and $N=2^{12},2^{16},\ldots,2^{64}$.
These longitudinal paths cover $n_T\geq n$. For $n_T<n$, we examine
selected allocations through fixed-source scans.
For these scans, we fix $n=2^{24},2^{40}$ and scan
$\kappa_T=\log n_T/\log n$ from $0.50$ to $1.75$, including the
polynomial and logarithmic branch boundaries.  Sample sizes are rounded
to integers.  Each fixed-$n$ scan uses one common numerical resolution,
so source-only risk is exactly invariant to $n_T$.
All 180 configurations use analytic expectations of the regression
estimators, including random-design variance and actual bias.
Raw-regression Monte Carlo checks only the six configurations
$n=n_T=2^{12},2^{16},2^{20}$ at $a=2,3$, with 400 repetitions each.
SAR intervals use 4000 paired bootstrap resamples, recomputing the library
maximum and channel minimum within each draw.  Larger-$N$ curves contain
no simulated losses or sampling intervals.

\begin{figure}[!t]
  \centering
  \includegraphics[width=\linewidth]{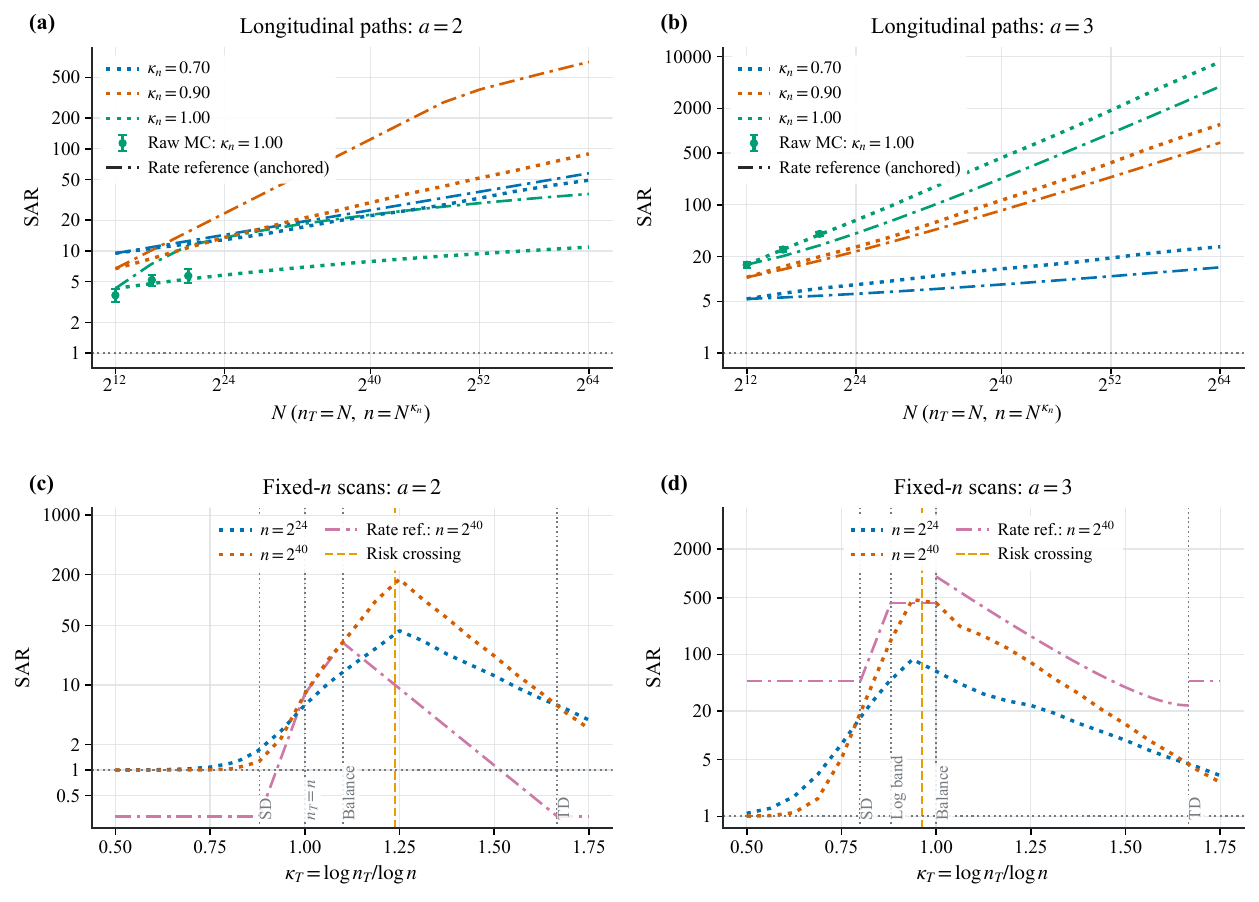}
  \caption{Integral-link SAR at $a=2,3$ with $\beta=1$ on the fixed
  three-function library.  Panels (a,b) use $n_T=N$,
  $n=N^{\kappa_n}$ for $\kappa_n=.70,.90,1.00$; dotted curves are
  analytic-risk SAR and same-color dash--dot curves are theoretical
  references aligned once at $N_0=2^{12}$.
  Panels (c,d) fix $n=2^{24},2^{40}$ and scan $\kappa_T$; purple
  dash--dot segments give the $n=2^{40}$ reference aligned once at
  $\kappa_T=1$.  Gray vertical lines mark theoretical boundaries and
  balances, while orange dashed lines mark analytic-risk channel
  crossings.  Markers at $\kappa_n=1$ use 400 raw-regression repetitions
  with pointwise 95\% percentile-bootstrap bars.  The horizontal line is
  SAR$=1$.}
  \label{fig:sim:iii-critical}
\end{figure}

Figure~\ref{fig:sim:iii-critical} reports the longitudinal paths at $a=2,3$
and fixed-source scans together with their theoretical rate references.

\paragraph{Theoretical SAR references.}
Let $s_a(n,n_T)=\min\{r_{S,a}(n),r_T(n_T)\}/r_{J,a}(n,n_T)$,
where $r_T=n_T^{-4/5}$, $r_{S,2}=(\log n)/n$,
$r_{S,3}=n^{-4/5}$, and $r_{J,a}$ is the complete unit-constant
piecewise expression in~\eqref{eq:iii-a2} or~\eqref{eq:iii-a3}
at $\beta=1$, including its logarithmic factors.
Each longitudinal reference is multiplied by a single constant to match
the corresponding analytic-risk SAR at $N_0=2^{12}$.
The scan reference uses $n=2^{40}$ and is aligned once at $\kappa_T=1$.
No horizontal fitting or separate branch alignment is applied.
At a branch interface we use the expression on the left for anchoring
and draw the one-sided reference segments separately: the order formulas
need not agree in their constants, particularly at $a=3$.
Their aligned constant levels outside the SLP window need not equal one.

\paragraph{Results at $a=2,3$.}
Figure~\ref{fig:sim:iii-critical} shows increasing longitudinal SAR and
hump-shaped allocation scans for these estimators.
The reference comparison reveals nonuniform rate agreement.
At $a=2$, the terminal analytic SAR divided by its anchored reference
is $.856$, $.126$, and $.302$ for $\kappa_n=.70,.90,1.00$;
the latter two paths grow substantially more slowly over this range.
At $a=3$ the corresponding ratios are $1.894$, $1.758$, and $2.139$,
so finite-range drift remains visible even when the broad growth pattern
agrees.
The scan-grid maxima occur at $\kappa_T=1.25$ for $a=2$ and approximately
$0.9375$ for $a=3$ at both fixed source sizes.
At $n=2^{40}$, interpolating the analytic log risk ratio
gives single-channel crossings at $1.238$ and $.962$, compared with
the unit-constant theoretical balances $1.100$ and $1.000$.
The differing risk constants can shift these crossings and the scan shape;
the references are not horizontally adjusted to remove that displacement.
At the right endpoint, SAR decreases from $3.89$ to $3.08$
for $a=2$ and from $3.18$ to $2.64$ for $a=3$ as the fixed source size
increases.

\paragraph{Mechanisms and normalization.}
The expected-risk decomposition identifies two distinct sources of the
logarithmic factors: accumulation in the compensated cumulative variance
at $a=2$, and accumulation across oscillatory shell blocks at $a=3$.
Appendix~\ref{app:sim:iii-normalization} reports these decompositions,
the remaining finite-sample corrections, and the normalized-risk checks.
Figure~\ref{fig:sim:iii-noncritical} concerns a fixed 64-dimensional span.
Figure~\ref{fig:sim:iii-critical} concerns specified estimators in a fixed infinite
orthonormal system, approximated at a controlled numerical resolution.

\FloatBarrier

\subsection{Unbounded derivative link: regression Monte Carlo and surrogate-risk diagnostics}
\label{app:sim:iv}

\begin{figure}[t]
  \centering
  \includegraphics[width=\linewidth]{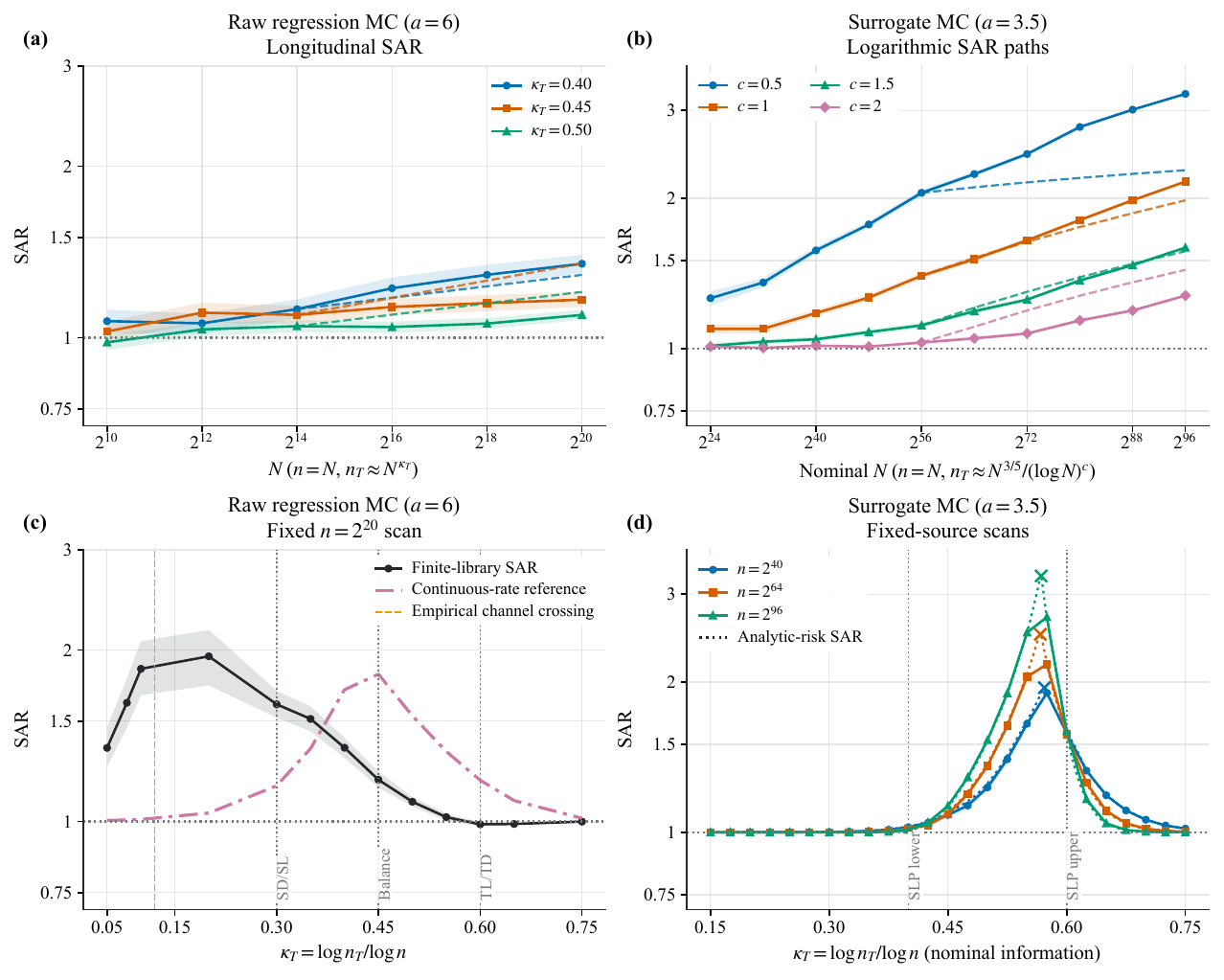}
  \caption{Derivative-link longitudinal SAR and full fixed-source scans with $\beta=2$.
  Columns distinguish raw-regression Monte Carlo ($a=6$, left)
  from surrogate-model Monte Carlo ($a=3.5$, right).
  (a) Finite-library SAR on $\kappa_T=.40,.45,.50$ inside the theoretical SLP
  window; dashed guides have slopes $1/30,1/20,1/30$.
  (b) Logarithmic SLP paths with $c=.5,1,1.5,2$ and varying effective
  $\kappa_T(N)$; dashed guides show continuous-rate SAR shapes.
  Same-color dashed guides in (a) and (b) are aligned at the first fitted scale.
  (c) The fixed-$n=2^{20}$ scan with a common truth library and source
  data, including the $n_T=2,3$ small-sample diagnostics.
  The purple dash--dot curve is the continuous local-rate
  SAR reference without horizontal calibration; orange dashed lines mark
  observed channel crossings.  Gray dotted vertical lines mark asymptotic
  balances and window endpoints.
  (d) Full-range scans at three fixed source-information levels; crosses mark
  analytic source/target equality and dotted lines mark $.4,.6$.
  Colored dotted curves in (d) give exact-risk SAR for the surrogate
  estimators on a finer scan grid, without sampling bands.
  Solid lines and markers summarize Monte Carlo estimates; bands are
  pointwise 95\% paired-bootstrap intervals within each model.
  The joint-estimator SAR is shown without truncation below one.
  In (b,d), nominal $N$ specifies the information level of the coefficient
  experiment.}
  \label{fig:sim:iv}
\end{figure}

\paragraph{Design.}
We set $\beta=2$, hence $s=1$ and $a_c=3.5$.  Source covariates have density
$h(x)=a x^{a-1}$ and target covariates are uniform.
Source responses observe $f$, whereas target responses observe $\dot f$.
The left column of Figure~\ref{fig:sim:iv}
reports the raw-sample experiment with
\[
 a=6,\quad n=N=2^{10},2^{12},\ldots,2^{20},\quad
 \kappa_T\in\{.40,.45,.50\}, \quad
 n_T=\max\{4,\operatorname{round}(N^{\kappa_T})\}.
\]
The longitudinal panel isolates the source-limited, equalizer, and target-limited
paths inside the theoretical SLP window $.30<\kappa_T<.60$.
The fixed-source panel retains the full sampled range, including the window
boundaries and both outer regimes.  Path selection follows the theoretical
classification, not the observed SAR or fitted slope.
All six $a=3$ paths, with $\kappa_T=.30,.45,.50,.55,.60,.75$, and the
$a=6$ boundary/outer paths $.15,.30,.60,.75$ are retained as
controls in Appendix~\ref{app:sim:iv-controls}.  We also evaluate a
fixed-source scan at $a=3$.
Each cell uses three distinct compact-bump scale stresses and a smooth
quadratic control, with a 1025-point integration grid.  For each estimator,
we first average its loss for each truth and then maximize over the common
four-truth library.  The SAR numerator is the minimum of the two resulting
single-channel risks; all these operations are repeated within each of
2000 paired bootstrap draws.  Horizontal raw scans fix $n=2^{20}$, the truth
library, and the source observations, and report individual risks as well as SAR.
The $a=6$ scan additionally includes $n_T=2,3$, corresponding to
$\kappa_T=.05$ and $\log_2(3)/20\approx.07925$, to expose the left side
of the finite-sample gain profile.  These two scan points are exempt from
the longitudinal minimum of four target observations. These points use
the same estimator to examine behavior with very small target samples.

At $a=3.5$, we simulate the separately defined diagonal Gaussian
coefficient model.  Its basis-derived information constants, common observations,
linear estimators, and restricted signal class are specified below.
At nominal $N=2^{24},2^{32},\ldots,2^{96}$, we follow
\[
 n=N,\qquad n_T\approx N^{3/5}/(\log N)^c,\qquad
 c\in\{.5,1,1.5,2\}.
\]
These positive-$c$ paths probe the logarithmic SLP mechanism.  Their effective
sample exponent is
\[
 \kappa_T(N)=\frac{\log n_T}{\log N}
 \approx\frac35-c\frac{\log\log N}{\log N};
\]
unlike the power paths at $a=6$, they do not hold $\kappa_T$ fixed.
The $c=0$ bounded-SAR control and the power controls
$n_T\approx N^{\kappa_T}$ with $\kappa_T=.30,.45,.70$ are displayed separately
in Appendix~\ref{app:sim:iv-controls}.
Fixed-source scans use $n=2^{40},2^{64},2^{96}$.  Their Monte Carlo grid has
spacing $.025$ in $\kappa_T$; analytic-risk calculations for this surrogate are
also evaluated on a fixed $.0025$ grid across the entire scan range.
No information constant is fitted to a desired peak.  Three fixed-source
sizes diagnose the slow finite-information displacement at $a=3.5$.

\paragraph{Raw-regression evidence and limitations.}
In Figure~\ref{fig:sim:iv}(c), the fixed-library source/target
crossing for $a=6$ is approximately $.120$, obtained by interpolating the
log risk ratio, rather than the asymptotic balance $.45$.
The $n_T=2,3,4$ points have SAR $1.347,1.616,1.854$,
revealing a rising left flank before the displayed grid maximum at
$n_T=16$.  The paired SAR ratio from $n_T=2$ to $3$ is
$1.200$ $[1.128,1.291]$.
The source risk is exactly constant within the scan.  Of the 78 longitudinal
raw cells, 77 are target-favoured; the sole source-favoured cell is the largest
$a=3,\kappa_T=.30$ configuration.  All $a=6$ longitudinal cells remain
target-favoured at the simulated sizes.

At $a=6$, $\kappa_T=.45$, terminal SAR is $1.166$ $[1.133,1.203]$.
Its tail-slope interval in Table~\ref{tab:sim:summary} excludes zero,
but the estimate lies below the asymptotic exponent $.05$ and the
finite-range continuous-rate reference. For the coefficient experiment,
the corresponding control at $a=6$
has slope $.0460$ $[.0448,.0472]$ and terminal SAR $7.99$.

\paragraph{Logarithmic paths at $a=3.5$ and remaining rate uncertainty.}
Figure~\ref{fig:sim:iv}(b) shows slowly increasing coefficient SAR on
the positive-$c$ paths; their empirical terminal SARs for
$c=.5,1,1.5,2$ are respectively $3.237,2.160,1.593,1.277$.
The asymptotic logarithmic equalizer $c=1.5$ is not the best of these
finite-range paths. At $n=2^{96}$, the surrogate's analytic source/target
risks cross near $\kappa_T=.568$.

Appendix~\ref{app:sim:iv-controls} reports the normalization checks and
remaining rate uncertainty at $a=3.5$ using the continuous reference $\rho_a$
defined below.  Independent grid-refinement diagnostics from 1025 to 2049 points
changes raw mean risks by at most $.0083\%$ on six 400-repeat diagnostic cells;
the coefficient omitted-energy bound is below $3.4\times10^{-8}$ of joint
risk, and a resolution-40 check finds no detectable floating-point change
in the analytic risks.

\subsection{Experimental details and numerical summary}
\label{app:sim:details}

\paragraph{Section~\ref{sec:slp4} stress library and mixed fit.}
Let $\phi(z)=(1-z^2)^4\mathbf 1\{|z|<1\}$.  The three scale stresses use
\[
 f_N(x)=\frac{b^2}{8}\sum_k(-1)^k
 \phi\!\left(\frac{x-c_k}{b}\right),\qquad
 \dot f_N(x)=\frac{b}{8}\sum_k(-1)^k
 \phi'\!\left(\frac{x-c_k}{b}\right).
\]
Here $b$ is the bump radius.  Centers start at $\delta L+b$, have spacing
$2.5b$, and end no later than $.92L-b$; the source-, target-, and mixed-scale
stresses have $\delta=.12,.20,.16$, respectively.  Since the supports are
disjoint and $\|\phi''\|_\infty=8$, each derivative is Lipschitz with
constant at most one.  Below $a_c$, the source scale is
$(L,b)=(1,.18n^{-1/5})$; above $a_c$, it is
$(L,b)=(n^{-1/(a+4)},.18n^{-1/(a+4)})$.
The target scale is $(1,.18n_T^{-1/3})$.  The mixed scale uses
$L=1$ below $a_c$ and $L=1\wedge(n_T^{5/3}/n)^{1/(a-1)}$ above it, with
$b=.18\min(n_T^{-1/3},L)$.  Thus a small boundary region is retained by
shrinking the radius, not by moving the support to the interior.
The fourth truth is $f(x)=x^2/4$.
For horizontal scans, all four truths are constructed once using
$n_T=\operatorname{round}(n^{.45})$ and held fixed.

The implemented degree-one local fits use an Epanechnikov kernel, exact
truncated-window masses, and a local bias--information bandwidth balance,
restricted numerically to $[10^{-4},.45]$.  Single-channel fits use the
corresponding source or target information.  A trace-scaled ridge of $10^{-9}$
stabilizes the normal equations.  Projection onto the known derivative bound
$[-1,1]$ is a statistical safeguard, not merely a numerical rounding operation.
All methods and truths share the same covariates and Gaussian noises within
each repetition.  In a fixed-source scan, the source observations and truth
library are identical at every $n_T$; target samples are nested.  Diagnostic
analyses separate conditional squared bias and noise risk before
projection from the reported projected risk.

\paragraph{Section~\ref{sec:slp4} surrogate model, Monte Carlo losses, and analytic risks.}
On shell $[L,2L]$, $L=2^{-j}$, use levels with cell width $b=L2^{-\ell}$
and multiplicity $m=2^\ell$, for $1\leq j\leq36$ and
$0\leq\ell\leq36-j$.  Set
$J_0=\int_{-1}^1\phi^2$, $J_1=\int_{-1}^1(\phi')^2$,
$K=J_0/(4J_1)=7/136$, $A=\sqrt{J_1}/(8\,2^{3/2})$, and
$\sigma=.25$.  Here $b$ denotes cell width, so its bump radius is $b/2$.
The normalized derivative mode has primitive energy $Kb^2$ and the
coefficient envelope $t_b=Ab^{3/2}$.  We explicitly replace the design
density by its shell average $\bar h_L=(2^a-1)L^{a-1}$ and use the diagonal
Gaussian experiment
\[
 Z_S=\theta+I_S^{-1/2}\xi_S,\qquad
 Z_T=\theta+I_T^{-1/2}\xi_T,\qquad
 I_S=\frac{n\bar h_LKb^2}{\sigma^2},\quad
 I_T=\frac{n_T}{\sigma^2},
\]
with independent standard normal coordinates. Shell averaging, diagonal
Gaussian observations, and the orthogonal coefficient loss define this
surrogate experiment. Equivalence to the
random-design regression model is not assumed.

All three implemented linear shrinkage methods use the same coordinates and
observations, with envelope-based factor $w(I)=t_b^2I/(1+t_b^2I)$.
Single-channel methods use their respective information; the joint method
uses total information and the information-weighted observations.
The sequence signal class activates at most one level per shell, with
$|\theta_k|\leq t_b$ on that level.  For each fixed estimator, choose on
every shell the level maximizing $mt_b^2(1-w)^2$, before generating data.
The resulting three truths attain the respective worst risks for these
estimators on this restricted class; every estimator is evaluated on all
three truths, and variance from inactive levels is included. Monte Carlo losses are
computed from common Gaussian signal-direction projections and their
orthogonal-noise Wishart sums, an exact aggregation of squared coefficient
errors.  They are not random draws from a prescribed theoretical risk.
We separately evaluate the same estimators' exact expected losses for
analytic-risk comparisons.
The unresolved signal-energy bound is $4A^2 2^{-72}$; extending the
coefficient resolution from 36 to 40 is checked separately.

\paragraph{Continuous rate reference.}
For Section~\ref{sec:slp4} normalization, solve the unrestricted exact local-mass
balance for $b_*(x)$ and integrate
$\rho_a(n,n_T)=\int_0^1b_*^2(x)\,dx$, including the boundary region.
This continuous diagnostic avoids jumps caused by choosing unit constants
independently in the branches of an $\asymp$ formula. It serves as the
theoretical rate benchmark for normalization. The comparison at $a=3.5$
that omits the logarithmic factor replaces it by $n^{-2/5}$. In all-path diagnostic
comparisons, both normalizations are anchored to one at the first scale,
without fitting slopes or channel-information constants.  Representative-path
comparisons use unrescaled $\widehat R_J/\rho_a$, matching the display convention
of Section~\ref{sec:slp3}.

\begin{table}[H]
\centering
\caption{Numerical summary of selected regression and coefficient-surrogate
rate paths.
We report the positive risk-decay exponent
$\widehat\alpha_R=-d\log(\widehat R)/d\log N$ and the SAR exponent
$\widehat\alpha_S=d\log(\widehat{\mathrm{SAR}})/d\log N$; brackets are
paired-bootstrap 95\% intervals.  Section~\ref{sec:slp2} integral rows use candidate
mean-risk minima with $\beta=1$ and the six scales $2^{15},\ldots,2^{20}$;
its derivative row uses $\beta=2$, raw joint risk, and
$2^{14},2^{16},2^{18},2^{20}$.  Section~\ref{sec:slp4} uses $\beta=2$.
Section~\ref{sec:slp4} fits use the last four raw
scales or six coefficient scales. Their reference slopes are computed from
the continuous local-rate integral and its one-channel analogues.  Other references are asymptotic exponents.}
\label{tab:sim:summary}
\renewcommand{\arraystretch}{1.12}
\scriptsize
\resizebox{\linewidth}{!}{%
\begin{tabular}{@{}llccccc@{}}
\toprule
Section & Path & $\widehat\alpha_R$ (95\% CI) & Reference
& $\widehat\alpha_S$ (95\% CI) & Reference & Terminal SAR \\
\midrule
\ref{sec:slp2} & $\kappa_n=.70$
& $.881\ [.861,.901]$ & $.850$
& $.041\ [.023,.057]$ & $.050$ & $1.78$ \\
\ref{sec:slp2} & $\kappa_n=.80$ (equalizer)
& $.930\ [.903,.958]$ & $.900$
& $.112\ [.065,.159]$ & $.100$ & $2.71$ \\
\ref{sec:slp2} & $\kappa_n=.90$
& $.966\ [.929,.997]$ & $.950$
& $.051\ [.009,.088]$ & $.050$ & $1.29$ \\
\ref{sec:slp2} & derivative control
& $.414\ [.402,.426]$ & $.400$
& $-.008\ [-.012,-.004]$ & $0$ & $1.00$ \\
\midrule
\ref{sec:slp4} & $a=3<a_c$, $\kappa_T=.60$
& $.412\ [.399,.424]$ & $.406$
& $.0012\ [-.0027,.0048]$ & $-.0019$ & $.99$ \\
\ref{sec:slp4} & $a=a_c$, $c=1.5$
& $.3920\ [.3917,.3922]$ & $.3962$ (log)
& $.0130\ [.0126,.0134]$ & $.0123$ (log) & $1.59$ \\
\ref{sec:slp4} & $a=6>a_c$, $\kappa_T=.45$
& $.331\ [.315,.345]$ & $.337$
& $.0145\ [.0032,.0261]$ & $.0280$ & $1.17$ \\

\bottomrule
\end{tabular}}
\end{table}

\subsection{Additional integral-link risk checks}
\label{app:sim:iii-normalization}
The normalization checks at $a=1.5,2.5,4$ divide mixed risk by the complete
Section~\ref{sec:slp3} rate. At the theoretical equalizers
$\kappa_n=0.8,0.9,1.2$ for $a=1.5,2.5,4$, the max/min ratios over the five
fitted tail scales are $1.10,1.15,1.18$ and the normalized-risk slopes are
$-0.039,-0.049,0.032$. The $a=2.5$ slope interval is $[-0.088,-0.009]$,
indicating residual downward drift.

\paragraph{Methods and risk calculation at $a=2,3$.}
The analysis at $a=2,3$ uses roots $1$ and $\sqrt3(2x-1)$ and two orthonormal
piecewise-linear wavelets on each dyadic interval, each with zero zeroth and
first moments.  Their nested piecewise-polynomial spaces make integrated
squared error additive without requiring independent coefficient errors.
The implemented methods use Taylor-compensated IPW and integration-by-parts
source statistics, direct target projection, and deterministic combination,
shrinkage, and cutoff rules based on variance and squared-bias bounds.
The bounds use conditional second moments at most one and the squared
coefficient envelope $(H/8)^2b^5$ at width $b$, rather than truth-specific
risks; the joint bias bound includes the source weight.

At $a=3$, the unregularized boundary wavelets have finite IPW variance
but an infinite fourth moment.  We therefore apply zeroth-order Taylor
compensation below $\delta=(H^2n)^{-1/5}$, capped at half the wavelet's
support width.  This rule was fixed before production Monte Carlo.
Actual signed boundary bias and modified variance enter the expected risk,
while combination and shrinkage use variance plus squared-bias bounds.
The bounded weights of each finite-resolution estimator yield finite
variance of the sampled losses.

Expected risks include signal-induced design variance, observation noise,
boundary and shrinkage bias, and the exact projection tail of the quadratic
truths.  Numerical resolution grows with the information range to
approximate one fixed infinite system, rather than changing the truth class.
The bound on the omitted tree contribution and projection tail is at most
$0.453\%$ of the corresponding fixed-library risk over the 180 cells.
Increasing resolution by two levels at seven prescribed configurations
changes the actual risks by at most $0.00478\%$.
The effective block count in the accompanying diagnostics is the participation count
$(\sum_d R_d)^2/\sum_d R_d^2$, not a count based on a fitted threshold.
The blocks combine orthogonal modes within physical shells and boundary
modes assigned by support scale.  Their physical supports need not be
disjoint, and their sampling errors need not be independent.

The analytic risk decomposition for these estimators separates cumulative
variance at $a=2$ from oscillatory shell accumulation at $a=3$. For the first cumulative functional
with kernel $K(t)=1-t$, the compensated source weight gives noise variance
$\sigma^2\{\tfrac12\log(1/b)+\tfrac14\}/n$.
At $a=3$ and $f=0$, each effective orthogonal shell block contributes
approximately $0.0127/N$ to combined oscillatory risk, and the participation
count increases from $5.89$ at $N=2^{20}$ to $13.85$ at $N=2^{60}$.
At $N=2^{64}$, the corresponding log-normalized values remain $26.7\%$
and $10.0\%$ above the predicted leading coefficients for $a=2,3$.

\paragraph{Monte Carlo checks and uncertainty.}
Each repetition generates independent source and target designs and
responses, applies the finite-resolution estimators, and computes squared
reconstruction error. Common random numbers pair the truths and methods.
The predeclared affine truth $f(x)=0.2+0.4x$ avoids a library-maximization
effect in direct mean-risk comparisons. At $a=3$, $N=2^{12}$, its
400-repeat source-only mean is $0.966$ times the analytic expectation.
The first cumulative cutoff contains only $0.074$ expected source
observations per repetition, making rare boundary observations influential.
A separate independent 400-repeat diagnostic gives ratio $0.975$ with relative
MC standard error $0.214$.

Agreement is not uniform across total-risk checks.
At $a=3$, $N=2^{20}$, zero-truth combined total risk is $0.909$ times its
expectation, or $-3.14$ estimated MC standard errors away. The corresponding
affine-truth ratio is $0.963$, or $-1.16$ estimated standard errors away.
At $a=2$, $N=2^{12}$, the selected cumulative component is $-2.84$ estimated
standard errors below its expectation.
The reported errors describe approximate pointwise MC precision, and rare
boundary contributions can make normal approximations inaccurate with 400
repetitions.

\subsection{Derivative-link negative controls and remaining rate uncertainty}
\label{app:sim:iv-controls}
\begin{figure}[t]
  \centering
  \includegraphics[width=\linewidth]{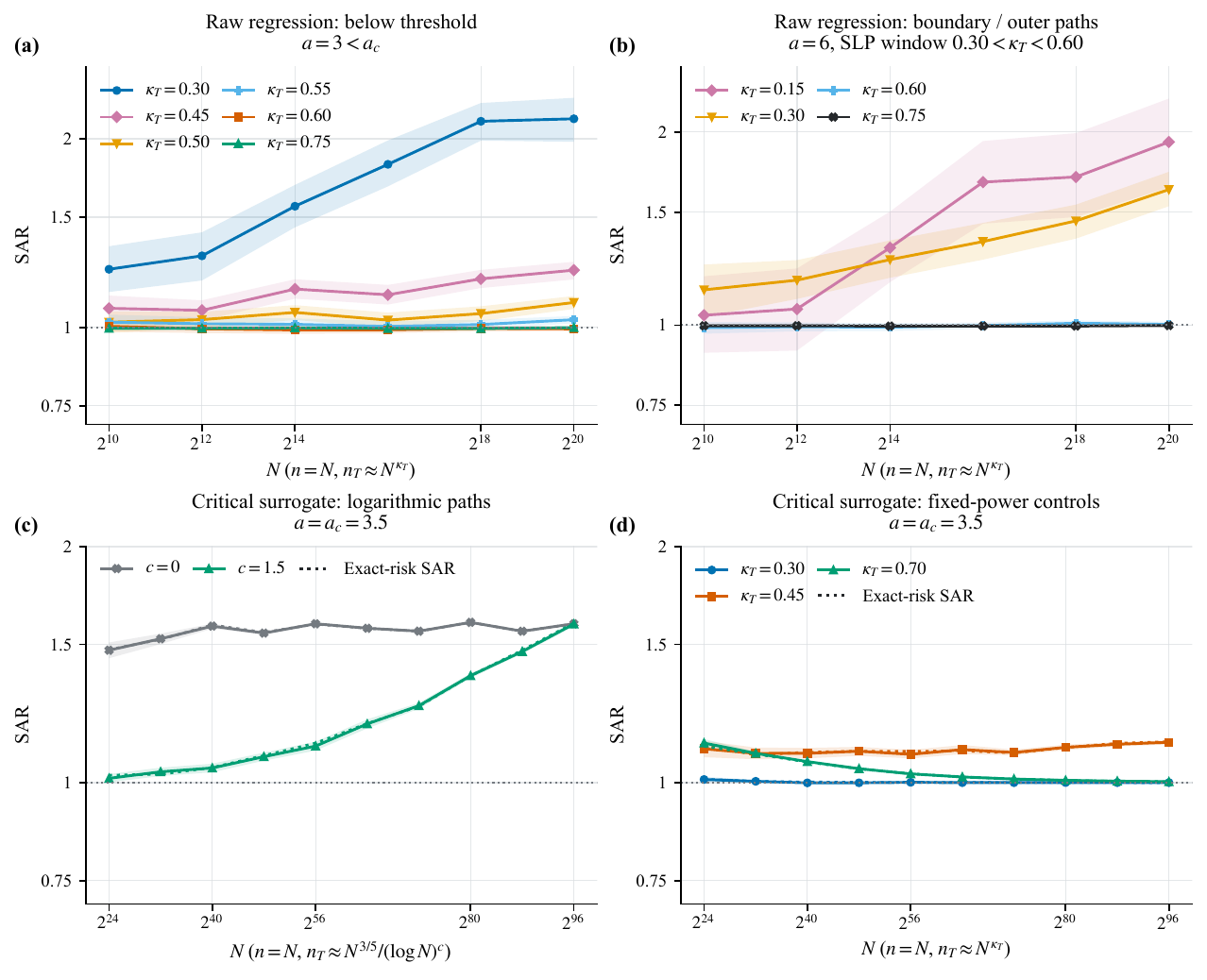}
  \caption{Derivative-link negative controls with $\beta=2$ and 400
  Monte Carlo repetitions per configuration.
  (a) All six raw-regression paths below the threshold, $a=3<a_c=3.5$.
  (b) All four raw-regression paths at or outside the $a=6$ SLP window
  $.30<\kappa_T<.60$.
  (c) Coefficient-surrogate logarithmic paths at $a=3.5$ with $c=0$ and $c=1.5$,
  where $n_T\approx N^{3/5}/(\log N)^c$.
  (d) All three fixed-power controls at $a=3.5$, $n_T\approx N^{\kappa_T}$.
  Solid curves and markers give MC SAR estimates with
  pointwise 95\% paired-bootstrap bands. Colored dotted curves in
  (c,d) are exact-risk SARs for the specified surrogate estimators, without
  sampling bands. Gray horizontal lines mark SAR$=1$.
  Below-threshold finite-range growth and SAR below one are shown without
  truncation. In (c,d), nominal $N$ denotes coefficient information.}
  \label{fig:sim:iv-negative-controls}
\end{figure}

Figure~\ref{fig:sim:iv-negative-controls}(a,b) shows the six $a=3<a_c$
raw-regression paths and four $a=6$ boundary and outer paths.
Below the threshold, the $\kappa_T=.55,.60,.75$ tail SAR-slope intervals
contain zero, whereas $.30,.45,.50$ have positive finite-range slopes
$.0809$ $[.0568,.1060]$, $.0193$ $[.0083,.0299]$, and
$.0096$ $[.0015,.0180]$. Three below-threshold paths thus show positive SAR
growth over the fitted range.
The complete $a=3$ scan has an observed single-channel crossing near
$.315$, displaced from the asymptotic balance $.60$.

At $a=3.5$, panel (c) pairs the $c=0$ bounded-SAR control with
the theoretical logarithmic equalizer $c=1.5$. Their terminal SARs are
$1.594$ and $1.593$, respectively. The other positive-$c$ paths remain in
Figure~\ref{fig:sim:iv}(b). Panel (d) retains all three fixed-power controls
$\kappa_T=.30,.45,.70$, which have no divergent SAR under the $a=3.5$
rate formula.

The remaining normalization checks use the full-domain reference $\rho_a$.
Over the last six scales at $a=3.5$, the max/min ratios of
$\widehat R_J/\rho_a$ for $c=0,.5,1,1.5,2$ are
$1.050,1.028,1.049,1.125,1.277$.
Without the logarithmic factor they are $1.050,1.095,1.143,1.244,1.426$.
The logarithmic reference reduces drift, but the normalized-slope intervals
for $c=1.5,2$ exclude zero.
\FloatBarrier

\FloatBarrier
\section{Empirical Illustrations with IO-VNBD}
\label{app:empirical}

The IO-VNBD benchmark contains synchronized 10-Hz vehicle CAN and smartphone
measurements collected on public-road trajectories \citep{onyekpe2021iovnbd}.
Yaw rate $\omega(t)$ and unwrapped heading $\theta(t)$ are approximately
related by
\[
  \theta(t)-\theta(0)\approx \int_0^t\omega(s)\,ds,
  \qquad \omega(t)\approx \theta'(t).
\]
We first use observed signals to illustrate operator asymmetry, then
calibrate a separate semi-synthetic benchmark under unbounded covariate shift.

\subsection{Real-Data Operator Asymmetry}
\label{app:real-data}

\paragraph{Data and calibration.}
The fully observed analysis uses vehicle CAN yaw rate and smartphone GPS
heading. After exact-payload deduplication and a prespecified link audit,
we retain 21 nonoverlapping 180-second windows from 11 signal-distinct track
files. The first 36 seconds form a calibration block and the remaining 144
seconds contain 1441 aligned analysis observations. Window selection uses
both calibration and analysis-block link diagnostics.

On the calibration block, let
$I(t)=\int_0^t\omega(s)\,ds$ and
$\Theta(t)=\theta(t)-\theta(0)$.  After centering both trajectories, we set
\[
 \widehat c=
 \frac{\sum_t\{I(t)-\bar I\}\{\Theta(t)-\bar\Theta\}}
      {\sum_t\{I(t)-\bar I\}^2}.
\]
This scale, including its sign, is fixed before the analysis block.  Writing
$T=144$ seconds and $x=(t-t_0)/T$, the forward task uses
\[
 f(x)=\frac{T\widehat c}{\pi}\omega(t),
 \qquad
 g(x)=\frac{\theta(t)-\theta(t_0)}{\pi},
 \qquad g(x)\approx\int_0^x f(u)\,du.
\]
The reverse task uses the same normalized signals with source level $g$ and
target derivative $f$.  Thus no synthetic responses or artificial response
noise are introduced.

\paragraph{Sampling and estimators.}
For integration we take $n=\operatorname{round}(N^\kappa)$, $n_T=N$, with
$\kappa\in\{0.5,0.8,1\}$; the theory-relevant middle path is $\kappa=0.8$.
For differentiation we take $n=N$ and
$n_T=\operatorname{round}(N^{0.6})$. In both directions
$N\in\{48,72,108,162,243,360,480\}$ and risk is held-out target MSE.
Within each analysis block, five interlaced folds are formed from consecutive
five-second dependence blocks. Repetition $r$ holds out fold $r\bmod5$ and
samples disjoint source and target training timestamps by temporal strata.
This preserves the quasi-regular recorded time design without using test labels.

For integration, target-only estimation is boundary-truncated local linear
regression with Epanechnikov kernel and bandwidth $1.1n_T^{-1/5}$.  The source
primitive integrates the piecewise-linear interpolant through the sampled yaw
labels.  The mixed estimator smooths the target residual between that
primitive and the observed heading, using correction bandwidth
$1\wedge(n/n_T)^{1/2}$.  The three-path sensitivity check uses this same raw
mixed estimator for every $\kappa$, rather than changing estimators at the
theoretical branch boundaries.

For differentiation, the operator-aware degree-one local fit combines source
level rows $[1,z]$ with target derivative rows $[0,b^{-1}]$.  At evaluation
location $x$, $b$ numerically solves
\[
 b^2\{nH(x,b)b^2+n_TH_T(x,b)\}=1,
\]
with the exact mass of the boundary-truncated uniform window.  A small ridge
$10^{-9}$ times the local Gram scale and clipping to $[-40,40]$ are used only
for numerical stability.  Forty percent of the sampled target labels form an
internal validation subset.  Starting from source-only, the selector changes
to target-only or joint only if the paired validation improvement exceeds one
estimated standard error; after selection, the chosen candidate is refitted
on all training labels.  The held-out fold is never used for this decision.

\begin{figure}[ht]
  \centering
  \includegraphics[width=0.94\linewidth]{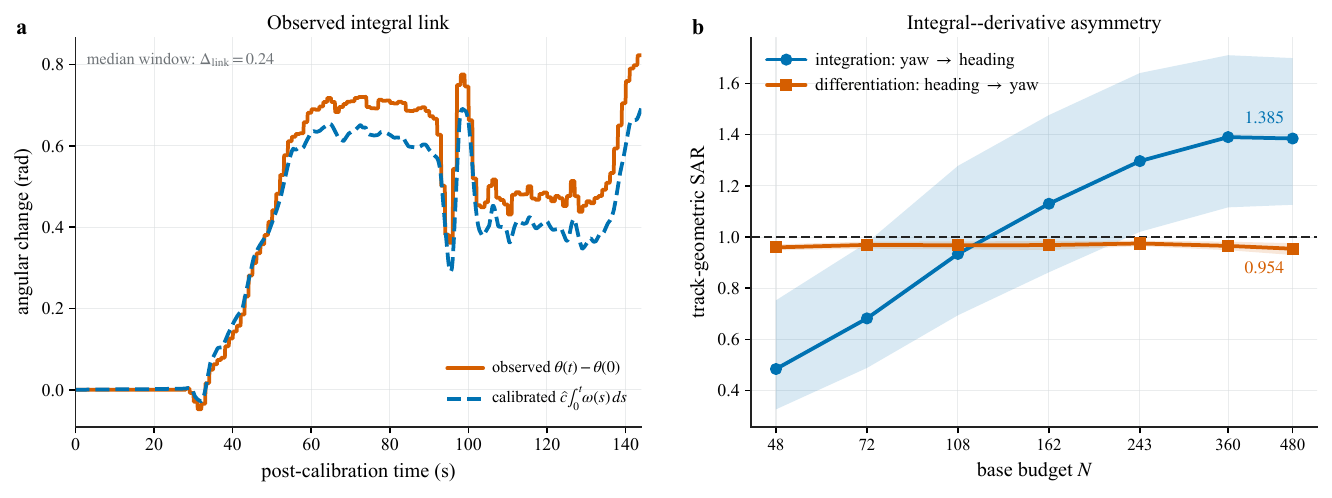}
  \caption{Real-data integral--derivative comparison on IO-VNBD.
  (a) Observed heading change and the calibrated yaw-rate integral in the
  deterministically selected median-discrepancy window.
  (b) Equal-track geometric SAR on the integration middle path and the
  safeguarded differentiation path. Bands are exploratory 95\%
  track-resampling intervals conditional on the audited windows.}
  \label{fig:real-asymmetry}
\end{figure}

\paragraph{Risk summaries and results.}
Each window--budget cell has 400 repetitions. We average first within window
and then within track, giving the 11 tracks equal weight. Pointwise 95\%
track-resampling intervals use 4000 bootstrap resamples of complete track-level
paths, preserving cross-budget and cross-estimator dependence. These exploratory
intervals condition on the audited windows and exclude window-selection
variability. The analysis assesses the qualitative integral--derivative
comparison in the observed signals.
Figure~\ref{fig:real-asymmetry} uses the track-geometric quantity
$\exp\{11^{-1}\sum_{j=1}^{11}\log(\mathrm{SAR}_j)\}$ for the operator
comparison.

As a descriptive sensitivity check, the middle integration path has the
largest terminal ratio-of-mean-risks SAR, $1.177$, whereas the $\kappa=0.5$
and $\kappa=1$ values are $1.064$ and $0.998$; all three pooled-risk
intervals include one. On the $\kappa=0.8$ integration path in
Figure~\ref{fig:real-asymmetry}(b), the terminal track-geometric SAR is $1.385$
with 95\% interval $[1.126,1.699]$, and 9 of 11 tracks favor mixing.  In the
differentiation direction, the corresponding safeguarded SAR is $0.954$
$[0.930,0.976]$, with only 1 of 11 tracks favoring transfer.  Eight tracks lie
in the joint quadrant ``integration gain, differentiation no gain.'' This
contrast is consistent with the operator asymmetry in
Remark~\ref{rem:asymmetry}.

\subsection{Semi-Synthetic Calibration under Unbounded Covariate Shift}
\label{app:semi-synthetic-unbounded}

The observed time design cannot isolate an unbounded density ratio. The
semi-synthetic benchmark therefore uses constructed covariates and latent
functions, with observation errors calibrated from vehicle GPS heading and
CAN yaw rate.

\paragraph{Residual calibration.}
We form a separate residual library from audited 600-second windows.  On the
experiment portion of each window, a cubic spline with nine interior knots is
fit to normalized heading.  Heading minus this spline gives the level-residual
pool, and calibrated yaw rate minus the spline derivative gives the
derivative-residual pool.  Each pool is centered and standardized to unit
root-mean-square magnitude within window.  Semi-synthetic responses add
independent draws from the relevant pool, scaled by $0.25$. Within each
window, the resulting errors retain the empirical marginal shapes and are
independent across observations and channels, with covariate-independent
variance.

\paragraph{Design.}
For integration and for the primary differentiation path, source covariates
follow $Q_S=\operatorname{Beta}(6,1)$ and target covariates are uniform, so
$\frac{dQ_T}{dQ_S}(x)=\frac{1}{6x^5}$
is unbounded as $x\downarrow0$.  We take
$N=2^{10},\ldots,2^{18}$ source observations and
$n_T=\operatorname{round}(N^q)$ target observations.  Risks are averaged over
repetitions and windows within track before the signal-distinct tracks receive
equal weight; the bootstrap preserves each complete sample-size path.

For integration, the source observes $g'$ and the target observes $g$.  We use
$q\in\{0.45,0.60\}$ and a seven-truth library containing the zero truth,
signed cumulative and interior alternatives, and their signed sums.
The implemented source-only, target-only, and mixed operator-aware ridge
methods share a finite dictionary within each configuration, evaluated through
derivatives for source data and function values for target data.  They use
squared-residual weights $0.25^{-2}$ and a squared coefficient-norm penalty
with weight one; single-channel fits omit the other channel. The rate
guarantee for the hierarchical estimator does not directly cover these
finite-dictionary ridge fits.

The reported risk is the maximum mean integrated
squared error over this library, and the empirical SAR is
$\min(R_S^{\mathcal L},R_T^{\mathcal L})/R_C^{\mathcal L}$.  This panel uses
23 windows from 15 signal-distinct tracks, 16 residual draws per window and
sample size, and 9999 paired track-bootstrap draws.

For differentiation, the source observes a level function $f$ and the target
observes $f'$.  This semi-synthetic panel uses a single-vertex
scale-dependent Assouad truth family, rather than the four-truth regression
library in Section~\ref{app:sim:iv}.
At each sample size, risk and SAR refer to that single truth.
The operator-aware mixed
local-polynomial estimator is compared with its two single-channel versions
on the primary path $(a,q)=(6,0.45)$.  The controls are the target-dominating
path $(6,0.60)$, the below-threshold path $(3,0.60)$, and the bounded-ratio
path $(1,0.60)$.  This panel uses 33 windows from 20 signal-distinct tracks,
20 residual draws per window and sample size, and 4000 paired track-bootstrap
draws.

The comparisons below concern finite-range performance of the specified
estimators, conditional on the residual windows. We analyze each direction
separately because the truth families, estimators, repetition counts, and
fitted slope ranges differ across panels.

\begin{figure}[tbp]
  \centering
  \includegraphics[width=0.94\linewidth]{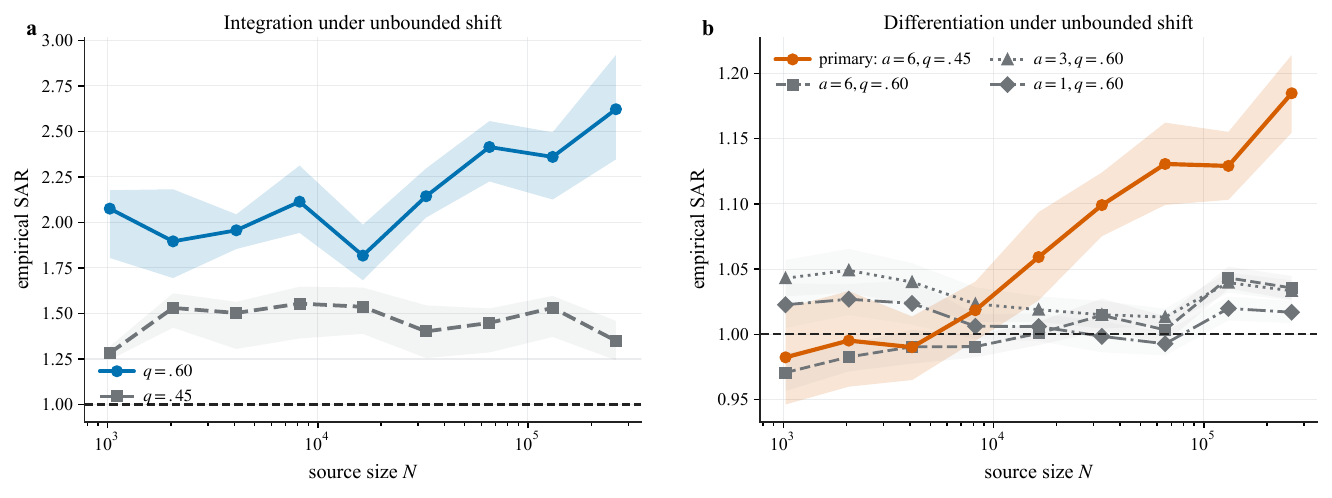}
  \caption{Semi-synthetic unbounded-shift calibration using empirical errors
  derived from the IO-VNBD signals introduced above.  (a) Integration SAR for
  the $q=.45$ and $q=.60$ paths.  (b) Differentiation SAR on the primary path
  and three controls for the sample-size regime, singularity threshold, and
  bounded density ratio.  Bands are pointwise 95\% signal-distinct-track
  bootstrap intervals, and the horizontal dashed lines mark SAR$=1$.}
  \label{fig:semi:unbounded-real}
\end{figure}

\paragraph{Results.}
Figure~\ref{fig:semi:unbounded-real}(a) shows the integration SAR paths.
On the $q=0.60$ path, the all-$N$ log-SAR slope is $0.0487$ with 95\%
interval $[0.0309,0.0665]$, and the terminal SAR is $2.62$
$[2.35,2.92]$.  The $q=0.45$ path has slope $0.0007$
$[-0.0116,0.0163]$ and terminal SAR $1.35$ $[1.24,1.46]$; the paired
$q=0.60$ minus $q=0.45$ slope contrast is $0.0480$
$[0.0220,0.0708]$.

Figure~\ref{fig:semi:unbounded-real}(b) shows the differentiation calibration.
On the primary $(a,q)=(6,0.45)$ path, the tail SAR slope is
$0.0323$ $[0.0161,0.0473]$, and the terminal SAR is $1.185$
$[1.153,1.213]$.  The primary slope exceeds the largest control slope by
$0.0178$ $[0.0001,0.0326]$. The contrast interval extends nearly to zero,
and the control paths also have small positive slopes over the observed range.

\end{document}